\input amstex.tex
\input amsppt.sty   
\magnification 1200
\vsize = 8.5 true in
\hsize=6.2 true in
\NoRunningHeads
\def\rightwave{{\ssize \sim\!\kern-.5pt\sim\!\kern-.5pt\sim\!\kern-.5pt\sim\!\kern-.5pt\sim\!\kern.5pt\sim
\!\kern-.5pt\sim\!\kern-.5pt\sim\!\kern-.5pt\sim\!\kern-.5pt\sim\!\kern-.5pt\sim\!\kern.5pt\sim}}
\def\tk{\tilde{\kern 1 pt\topsmash k}}
\def\bl{\bar{\kern 1.6pt\topsmash\ell}}
\def\bk{\bar{\kern 1.5pt\topsmash k}}
\def\td{\tilde{\kern -3pt \topsmash d}}
\def\bd{\bar{\kern -3pt \topsmash d}}

\def\ve{\varepsilon}

\def\supp{\text {supp\,}}
\def\dist{\text {dist\,}}
\def\mes{\text {mes\,}}
\def\ds{\raise 4pt\hbox{$\to{}$}\!\!\!\!\!\!\!\!ds}

\def\1{1 \!\!\!1}

\def\k{\lower 3pt\hbox{$\bar$}k}

\def\kk{\overset{\raise3pt\hbox{$\scriptscriptstyle\approx$}}\to {\topsmash k}}

\font\biggermath=cmex10 scaled\magstep4
\def\frownover#1{\lower4pt\hbox{$\topshave{\operatornamewithlimits{#1}
\limits^{\hbox{\biggermath\char'147}}}$}}

\TagsOnRight
\NoBlackBoxes
\parskip=\medskipamount
 at 10truept

\topmatter
\title 
Anderson Localization for Time Quasi Periodic \\
Random Sch\"odinger and Wave Operators\endtitle
\author{Jean~Bourgain and Wei-Min~Wang}
\endauthor
\address{Institute for Advanced Study, Einstein Drive,
Princeton, N.J. 08540, U.S.A.}
\endaddress
\email
{bourgain\@math.ias.edu}
\endemail
\address
{Department of Mathematics, Princeton University and
Institute for Advanced Study, Einstein Drive,
Princeton N.J. 08540, U.S.A.}
\endaddress
\email
{wmwang\@math.ias.edu}
\endemail
\thanks  Wei-Min Wang thanks A. Soffer and T. Spencer for many useful discussions and for initiating her to the subject.
The support of NSF grant DMS 9729992 is gratefully acknowleged. 
\endthanks 

\abstract
We prove that at large disorder, with large probability and for a set of
Diophantine frequencies of large measure, Anderson localization in
$\Bbb Z^d$ is {\it stable} under localized time-quasi-periodic perturbations by
proving that
the associated quasi-energy operator has pure point spectrum. The main tools are
the Fr\"ohlich-Spencer mechanism for the random component and
the Bourgain-Goldstein-Schlag mechanism for the quasi-periodic component. The
formulation of
this problem is motivated by questions of Anderson localization for
non-linear
Schr\"odinger equations.
\endabstract
\endtopmatter
\bigskip

\item {1.} Introduction
                                                                                                   
\item {2.} Exponential decay of the Green's functions of Schr\"odinger operator at fixed $E$ and $x$

\item {3.} Exponential decay of the Green's functions of Schr\"odinger operator at fixed $E$ and
$\theta$

\item {4.} The elimination of $E$ and frequency estimates

\item {5.} Proof of Anderson localization for the Schr\"odinger operator

\item {6.} Proof of Anderson localization for the wave operator

\item {7.} Appendix: Localization results for random Schr\"odinger operators

\bigskip

\heading I. Introduction \endheading

We prove persistence of Anderson localization for random Schr\"odinger and
random
wave operators under localized time-quasi-periodic perturbations. Given an initially
localized wave packet, Anderson localization is, roughly speaking, the
phenomenon
that the wave packet remains localized for all time. Schr\"odinger equation
is the
following:
$$
i\frac{\partial}{\partial t}\psi=(\epsilon\Delta+V)\psi,\tag 1.1$$
on $\Bbb R^d\times[0,\infty)$ or $\Bbb Z^d\times[0,\infty)$,
where $\epsilon>0$ is a parameter, $\Delta$ is the Laplacian (continuum or
discrete), $V$ the potential is a multiplication operator; wave equation is
$$\frac{\partial^2}{\partial t^2}\psi=(\epsilon\Delta+V)\psi,\tag 1.2$$
on $\Bbb R^d\times[0,\infty)$ or $\Bbb Z^d\times[0,\infty)$,
where the right hand side (RHS) is the same as in (1.1). In this paper, we
consider $V$ random, to be defined shortly.

When $V$ is independent of time, Anderson localization reduces to
prove that the {\it time independent} Schr\"odinger operator:
$$H_0=\epsilon\Delta+V,\,\tag 1.3$$ on $L^2(\Bbb R^d)$ or $\ell^2({\Bbb Z}^d)$
has pure point spectrum with exponentially localized (or sufficiently fast
decaying)
eigenfunctions. $0<\epsilon\ll 1$ is the large disorder case.

Anderson localization for {\it time independent}
random Schr\"odinger (or wave operator) at large disorder has been well
known
since the seminal work of Fr\"ohlich-Spencer \cite{FS}. It is a topic with
an
extensive literature \cite{GMP, FMSS, vDK, AM, AFHS, AENSS}, to name a few.

(Time independent) quasi-periodic Schr\"odinger operators in one dimension
are now
 well understood following the works in \cite{BG, FSWi, J, Sa, Sin} and
the related
works \cite{HS1,2}.
Recently in their fundamental paper \cite {BGS}, Bourgain-Goldstein-Schlag
proved Anderson
localization in two dimensions at large disorder under appropriate
arithmetic conditions on the frequency vector. (See \cite{Bo} for an
excellent
review and also overview of the subject and related things.) The papers
\cite{BG, BGS} play a central role in the construction here.

Below we specialize
to discrete random Schr\"odinger operator. $H_0$ is then defined as the
operator:
$$H_0=\epsilon\Delta+ V,\, \text {on}\, \ell^2({\Bbb Z}^d),\tag 1.4$$
\noindent
where the matrix element $\Delta_{ij}$, for $i$, $j\in\Bbb Z^d$ verify
$$\aligned \Delta_{ij}&=1\quad |i-j|_{\ell^1}=1\\
&=0\quad \text{otherwise};\endaligned\tag 1.5$$
\noindent
the potential function $V$ is a diagonal matrix:
$V=\text{diag}(v_j),\, j\in{\Bbb Z}^d$, where $\{v_j\}$ is a family of
independently
identically distributed (iid) real random variables with distribution $g$.
From now on, we write $|\,|$ for the $\ell^1$ norm: $|\,|_{\ell^1}$ on
$\Bbb Z^d$.
We denote $\ell^2$ norms by $\Vert\,\Vert$.
The probability space $\Omega$ is taken to be $\Bbb R^{\Bbb Z^d}$ and the measure
$P$ is
$\prod_{j\in {\Bbb Z}^d}g(dv_j)$.

As is well known, $\sigma(\Delta)=[-2d, 2d]$. Let $\supp\, g$ be the support
of
$g$, we know further (see e.g., \cite{CFKS,PF}) that
$$\sigma(H)=[-2d, 2d]+\gamma\supp\, g\quad a.s.\tag 1.6$$

The basic result proven in the references mentioned earlier is that under
certain regularity
conditions on $g$, for $0<\epsilon\ll 1$, and in any dimension $d$, the
spectrum of $H_0$ is
almost surely pure point with exponentially localized eigenfunctions, i.e.,
Anderson
localization, after the physicist P. W. Anderson
\cite{An}. Physically this manifests as a lack of conductivity due to the
localization of electrons. Anderson was the first one to explain this
phenomenon on
theoretical physics ground.

The study of electron conduction is a many body problem. One needs to take
into account
the interactions among electrons. This is a hard problem. The operator $H_0$
defined
in (1.4) corresponds to the so called 1-body approximation, where the
interaction
is approximated by the potential $V$. The equation governing the system is
(1.1) on $\Bbb Z^d\times[0,\infty)$.

As an approximation to the many body problem, when the interaction among
electrons
are weak, one studies the following {\it non-linear} Schr\"odinger equation
(cf \cite{DS, FSWa}):
$$i\frac{\partial}{\partial
t}\psi=(\epsilon\Delta+V)\psi+\delta|\psi|^p\psi,\qquad
(0<\delta\ll 1,\, p>0)\tag 1.7$$ on $\Bbb R^d\times[0,\infty)$ or $\Bbb Z^d\times[0,\infty)$.
In \cite{AF, AFS}, solutions to the non-linear eigenvalue problem corresponding to (1.7)
were found, which could be used to construct {\it time periodic} solutions where higher harmonics are absent. (See 
\cite{BFG}, for a Nekhoroshev type theorem in a related classical setting.)  But in order to obtain {\it time quasi
periodic} solutions or more general time periodic solutions to (1.7), one needs to study the corresponding 
{\it time dependent} random Schr\"odinger operators. 
 
We
remark here that the
non-linear Schr\"odinger equation in (1.7) is distinct from other more
commonly studied
non-linear Schr\"odinger equations in that the linear equation itself
already has
small-divisor problems. When $p=2$, (1.7) is also called the Gross-Pitaevskii equation,
which arises in the theory of vortices in boson systems \cite{Gr, Pi}.

In \cite{SW}, time-periodic, spatially localized perturbations of random
Schr\"odinger
operators were considered. It is proven that Anderson localization is {\it
stable}
under such perturbations.

In this paper, we prove that Anderson localization is also {\it stable}
under time-quasi-periodic, spatially localized perturbations with large
probability
and for a set of Diophantine frequencies of large measure. The techniques
here are
more involved than that in \cite{SW} as one needs to take care of the small
divisor
problem coming from the random component and the quasi-periodic component
simultaneously.

To be precise, we study the following time-quasi-periodic random Schr\"odinger
equation:
$$i\frac{\partial}{\partial t}\psi=(\epsilon\Delta+V+\Cal W)\psi\tag 1.8$$
and the time-quasi-periodic random wave equation
$$\frac{\partial^2}{\partial t^2}\psi=(\epsilon\Delta+V+\Cal W)\psi\tag
1.9$$
on $\Bbb Z^d\times [0,\infty)$, where as in (1.4), $V=\{v_j\}$ is a family of
(time-independent) i.i.d. random variables;
$$\Cal W=\Cal W(t,j)=\sum_{k=1}^{\nu}W_k(j)\cos 2\pi(\omega_k
t+\theta_k),\tag 1.10$$
where $$\aligned \omega&=(\omega_1,\cdots,\omega_{\nu})\in (0,1]^{\nu}\\
 \theta&=(\theta_1,\cdots,\theta_{\nu})\in (0,1]^{\nu}.\endaligned$$

To proceed further, we assume
\item {(H1)} $W_k(j)$ is such that, $$\sum_{k=1}^{\nu}|W(j)|\leq 2\nu\delta
e^{-b|j|}\,
(0<\delta\ll 1,\, b>0)$$
\item {(H2)} $\omega$ satisfies a Diophatine condition, 
$$\vert|n\cdot\omega\vert|_{\Bbb T^\nu}\geq\frac{c}{|n|^A}\qquad (n\neq 0,\, c>0,\, A>0).$$
We write $\omega\in\text{DC}_{A,c}$.
\item {(H3)} The probability distribution $g$ has bounded support, without
loss we assume
$\supp g\subset[-1, 1]$
\item {(H4)} $g$ is absolutely continuous with a bounded density $\tilde g$:
$$g(dv)=\tilde g(v) dv, \, \Vert \tilde g\Vert _{\infty}<\infty$$

\noindent
{\it Remark.} Some fast enough polynomial decay for
$\sum_{k=1}^{\nu}|W(j)|$
suffices. We assume (H1) for ease in writing.

It is well known (see e.g., \cite{Be, Ho1, 2, JL, Ya}) that the study of 
time-quasi-periodic equations like (1.8, 1.9) can be reduced to the spectral study
of their corresponding quasi-energy operators; here 
$$
K=\sum_{k=1}^{\nu}i{\omega_k}\frac{\partial}{\partial\theta_k}+\epsilon\Delta+
V+\sum_{k=1}^{\nu}W_k(j)\cos 2\pi\theta_k\tag 1.11
$$
for the Schr\"odinger equation in (1.8); and
$$K_w=(\sum_{k=1}^{\nu}{\omega_k}\frac{\partial}{\partial\theta_k})^2+\epsilon\Delta+
V+\sum_{k=1}^{\nu}W_k(j)\cos 2\pi\theta_k\tag 1.12
$$
for the wave equation in (1.9), on $\ell^2(\Bbb Z^d)\times L^2(\Bbb T^{\nu})$, (cf.
\cite{SW}).

We prove that with large probability and for a set of Diophantine
frequencies
$\omega\in (0,1]^{\nu}$ with large measure, both $K$ and $K_w$ have pure
point
spectrum with exponentially decaying (in the $\Bbb Z^d$ direction)
eigenfunctions.
(For a precise statement, see the Theorem in sect. 5 on $H$ and $H_w$ which
are
unitary equivalents of $K$ and $K_w$ respectively.) This implies in
particular
(after some standard gymnastics) that with large probability, for a set of
Diophantine
frequencies of large measure, and initial conditions $\psi(0)$ which are
localized in
space, the time evolutions $\psi(t)$ of (1.8) and (1.9) are almost periodic,
a.e. $\theta$,
(cf. e.g., \cite{SW, JL}).

We spare a few lines on the proof of Anderson localization for the unitary equivalents
$H$, $H_w$, which are obtained from $K$, $K_w$ by a partial Fourier transform in $\theta\in\Bbb T^\nu$,
(see (2.2, 6.2) for the precise expressions). Let $n$ be the dual variable of $\theta$, $n\in\Bbb Z^\nu$.
We know that for $0<\epsilon\ll 1$, roughly speaking, the Green's function decays exponentially 
in the $j$ directions, $j\in\Bbb Z^d$, due to Anderson localization of the original unperturbed
operator $H_0$ defined in (1.4). To prove Anderson localization for the perturbed
operators $H$, $H_w$ on $\ell^2(\Bbb Z^{d+\nu})$, we also need to prove exponential decay in the $n$ directions
using quasi-periodicity. This is however the ``classical" picture, as the quasi-periodic perturbation
{\it does not} commute with $H_0$.

To prove Anderson localization for $H$, $H_w$ on $\ell^2(\Bbb Z^{d+\nu})$, we put
the small-divisor problems originating from the random and quasi-periodic
components on equal footing and deal with them concurrently. For
the random component, we use the Fr\"ohlich-Spencer (FS) approach. The version
which is
well adapted to our purpose is the one in \cite{vDK}, summarized in the
appendix.
For the quasi-periodic component, we rely on semi-algebraic considerations,
Cartan
type of theorems for analytic matrix valued functions developed in the
series of
papers by Bourgain, Goldstein and Schlag (BGS) \cite{BG, BGS}, (see also \cite{Bo}). (The dynamics
here is simpler than that in \cite{BGS} due to the special quasi-periodic structure
of $H$, $H_w$.) The Diophantine
frequencies which
are excluded result from a Melnikov type of non-resonant conditions, (see
Lemmas 2.3, 6.1
(2.26-2.28, 6.10)).

Finally, for the experts, we wish to add that the constructive aspect of the BGS
mechanism is a more robust version of the FS mechanism. In BGS, at each scale, the number
of resonant sub-regions of the previous scale can grow sub-linearly; while in FS, at 
each scale, the number of resonant sub-regions of the previous scale is fixed
(see\cite {vDK}). In the quasi-periodic setting, one typically falls into the BGS
scenario.
\bigskip

\bigskip
\heading
{2. Exponential Decay of the Green's function of Schr\"odinger
operator at fixed $E$ and $x$}
\endheading

Recall from sect. 1, the quasi-energy operator $K$:
$$
K=\sum^\nu_{k=1} \omega_k\frac 1i \frac \partial{\partial\theta_k} +\ve\Delta +V+\sum^\nu_{k=1} W_k\cos 2\pi\theta_k
\tag 2.1
$$
on $\ell^2(\Bbb Z^d)\times L^2(\Bbb T^\nu)$, where 
$\omega=(\omega_1, \omega_2\cdots  \omega_\nu)\in (0, 1]^\nu$.
$V$ is the random potential on $\Bbb Z^d$, $0<\epsilon\ll 1$ and $W_k$ satisfies the decay properties
specified in (H1).

Performing a partial Fourier series transform, in the $\Bbb T^\nu$ variables only, we are led to study the following unitarily
equivalent operator:
$$
H=\delta_j\tilde\Delta_n+ n\cdot \omega+\ve \Delta_j+V_j
\tag {2.2}
$$
on $\ell^2(\Bbb Z^{d+\nu})$, where

\itemitem {$\bullet$}  \  $n\in \Bbb Z^\nu, j\in\Bbb Z^d$

\itemitem
{$\bullet$}  $\delta_j\tilde\Delta_n \, \overset {\text {def}}\to = \, \sum^\nu_{k=1} W_k(j) \Delta_k$ is an operator
on $\ell^2(\Bbb Z^\nu), \Delta_k$ is the standard discrete Laplacian on the $k^{\text th}$ copy of $\Bbb Z$.
$$
\Vert\delta_j\tilde\Delta_n\Vert_{\ell^2(\Bbb Z^\nu)} \leq 2\nu\delta e^{-b|j|}\qquad (b>0)
\tag{2.3}
$$
by using (H1).
(Some fast enough polynomial decay suffices.
We assume (2.3) for ease in writing.)

\itemitem {$\bullet$}  $\ve \Delta_j+V_j \, \overset {\text {def}}\to  =\, \ve \Delta+V$, we put in the subscript 
$j$ to stress that it came from an operator on $\ell^2(\Bbb Z^d)$

\itemitem {$\bullet$}   For simplicity, we now drop the tilde on $\Delta_n:\tilde \Delta_n \, \overset {\text {def}}\to =
\, \Delta_n$.

\itemitem {$\bullet$}
\ For $\Lambda \subset\Bbb Z^{d+\nu}$, $H_\Lambda$ is the restriction of $H$ to $\Lambda$: 
$$
H_\Lambda (j, n; j', n')\overset\text {def}\to = 
\cases H(j, n, j', n') &{ \text { if }} \ j, n\in \Lambda\,{\text{and}}\,  j', n' \in \Lambda \\
0 \qquad &{ \text { otherwise }}. \endcases\tag 2.4
$$

Let $X\subset \Bbb R^{\Bbb Z^d}$ be the set where the random Schr\"odinger operator $H_j =\ve\Delta_j+V_j$ exhibits
Anderson localization in a sense to be made precise in (2.14) where theorem 2.2 of [vDK], restated here as Theorem A is
applicable.
We note here only that since we require finite scale information, 
$$
\mes \,X< 1,\quad  \mes\, X \gtrsim 1- \frac 1{L^a} \, {(a>0)}
$$ 
where $L$ is the initial scale.

\noindent
{\bf 2.1 The initial estimate $(0^{\text {th}}$ step)}

Fix an energy $E$, fix $x\in X$, so that $H_j$ has Anderson localization.
Let $\theta\in \Bbb R$ and define
$$ 
H(\theta)=\delta_j\Delta_n+(n\cdot\omega+\theta)+\ve \Delta_j+V_j
\tag {2.5}
$$
on $\ell^2 (\Bbb Z^{d+\nu})$.
We study the Green's function
$$
 G_{\Lambda_0} (\theta, E)= (H_{\Lambda_0}(\theta)-E)^{-1}
\tag {2.6}
$$
for some $\Lambda_0=[-N_0, N_0]^{d+\nu}$, where $N_0$ is to be determined.
We call $\Lambda_0$, an $N_0$-box.
We do perturbation theory and for the $0^{\text {th}}$ step, we drop $\delta_j\tilde\Delta_n$.
We have after diagonalization
$$
H_{\Lambda_0, 0} \ \overset {\text {def}}\to = \, n\cdot \omega+\theta +\mu_j
\tag {2.7}
$$
where $\mu_j$ are the eigenvalues of $H_j$.
Since $\Vert\delta _j \Delta_n\Vert_{\ell^2(\Bbb Z^\nu)} \leq 2\nu\delta e^{ -b|j|}
$, from {(2.3)} we require that
$$
|n\cdot\omega+\theta+\mu_j-E|> 2c\delta
\tag{2.8}
$$
for some $c>2\nu$ and all $(n, j)\in\Lambda$.
So we estimate
the measure of the set of $\theta$ such that
$$
|n\cdot\omega+\theta+\mu_j-E|\leq 2c\delta
\tag {2.9}
$$
for some $(n, j)\in \Lambda$.
In this particularly simple case, we obtain
$$
\align
\mes \, \{ \theta | \Vert(H_{\Lambda_0, 0} -E)^{-1}\Vert &\geq (c\delta)^{-1} \}
\leq 4c \delta|\Lambda_0|^2\\
&=4c \delta(2N_0+1)^{2(d+\nu)}\tag{2.10}
\endalign
$$
Let 
$$
\sigma\in (0, 1), N_0 =[|\log c \delta |^{1/\sigma}]+1 
\tag {2.11}
$$
 ($[\cdot]$ is the integer part)
and $\Cal B_x (\Lambda_0, E)$ be the set defined in the left hand side of {(2.10)}.

We note that 
$$
\mes \, \Cal B_x (\Lambda_0, E) \leq e^{-\frac {N_0^\sigma}2}
\tag {2.12}
$$
for $N_0$ satisfying {(2.11)} and $0<\delta \ll1$.

\proclaim
{Lemma 2.1} 
There exists $\gamma'> 0$, such that for $\delta\ll 1$, on $\Bbb R\backslash\Cal B_x(\Lambda_0, E)$
$$
\align
&\Vert G_{\Lambda_0}(\theta, E)\Vert < e^{N_0^\sigma} \qquad \sigma \in (0, 1)\\
&|G_{\Lambda_0}(\theta, E)(m, m')|< e^{-\gamma'|m-m'|}\tag 2.13
\endalign
$$
for all $m, m'\in \Lambda, |m-m'|> N_0/4$.
\endproclaim

\noindent
{\bf Proof.}
The first inequality of (2.13) is a restatement of  (2.10), (2.11) .
To obtain the second inequality, we use the conclusion of Theorem 2.2 of [vDK] restated here as Theorem A
for the scale $N' =[N_0^{1/\alpha}]+1 \ (1<\alpha<2)$.
Theorem A states that the set $X_{N_0}$ where there is only one pairwise disjoint bad $N'$ box contained in the $N_0$ box
$[-N_0, N_0]^d$ has measure:
$$
\mes\, X_{N_0}\geq 1-\frac {N_0^{2d}}{{N'}^{2p'}}
\geq 1-\frac 1{N_0^{(2p'/\alpha - 2d)}} \quad (p'>0, 1<\alpha< 2).
\tag 2.14
$$
Fix $x\in X_{N_0}$, (assuming $2p'/\alpha -2d\gg 1$), using the resolvent expansion and the first equation of (2.13)
for the bad $N'$ box, we obtain that $\exists \gamma'>0$, such that on $\Bbb R\backslash \Cal
B_x(\Lambda, E)$,
$$
|(H_{j, \Lambda}(\theta)-E)^{-1} (i, i') |\leq e^{-\gamma'|i-i'|}\tag 2.15
$$
for all $i, i'\in [-N_0, N_0]^d$ and $|i-i'|> N/4$.

The second equation of (2.13) follows from Neumann series (in the $n$-direction), {(2.15)} and the decay condition
on $\delta_j$ in (2.3).
$\hfill\square$

\noindent
{\bf 2.2 A Wegner estimate (in $\theta$) for all scales}

We now prove an apriori estimate on $\Vert (H_\Lambda(\theta)-E)^{-1}\Vert$ for all finite subsets $\Lambda \subset\Bbb
Z^{d+\nu}$.
This estimate uses the special structure of (2.5) and hence holds only for Schr\"odinger and not for wave equations
e.g.
For those more general situations, we need to resort to Cartan-type of theorem for analytic matrix valued functions 
a la [BGS].
(For the experts, this saves us one subroutine and moreover we only need to work with cubes in $\Bbb Z^{d+\nu}$.)
Wave equation will be treated in sect 7.

\proclaim
{Lemma 2.2}
Let $E\in I$, an interval of length $\Cal O(1)$.
Let $\Lambda $ be a finite set in $\Bbb Z^{d+\nu}$.
$$
{\text{\rm mes\,}} \{\theta|\text{\rm dist\,} (E, H_\Lambda(\theta)) \leq \kappa)\} \leq C|\Lambda|\kappa
\tag{2.16}
$$
\endproclaim

\noindent
{\bf Proof.}
Let $\Cal N(\theta, \lambda)$ be the $\#$ of eigenvalues of $H_\Lambda(\theta)\leq \lambda$
$$
\align
&\mes\, \{\theta|\dist (E, H_\Lambda(\theta)\leq \kappa)\}\\
&\leq\int (\Cal N(\theta, E+\kappa)- \Cal N(\theta, E-\kappa))d\theta\\
& =\int_{|\theta|\lesssim \Cal O(1)N}\big(\Cal N(\theta, E+\kappa)-\Cal N(\theta, E-\kappa)\big)d\theta,
\tag{2.17}
\endalign
$$
since $\Cal N(\theta, E+\kappa) =\Cal N(\theta, E-\kappa)$ for $|\theta|> \Cal O(1)N$.
In view of {(2.5)}
$$
\Cal N(\theta, E \pm \kappa)=\Cal N(\theta \mp \kappa, E).\tag {2.18}
$$
Substituting (2.18) into (2.17)
we obtain
$$
\align
{(2.17)}& =\int(\Cal N(\theta-\kappa, E)- \Cal N(\theta+\kappa, E)) 
d\theta\\
 &=\int^{\Cal O(1) N-\kappa}_{-\theta(1)N-\kappa}\Cal N(\theta, E)d
\theta -\int^{\Cal O(1)N+\kappa}_{-\Cal O(1)N+\kappa}
\Cal N(\theta, E) d\theta\\
&\leq C|\Cal N(\theta, E)|_\infty \cdot \kappa\\
&\leq C|\Lambda|\kappa
\tag{2.19}
\endalign
$$
where we used the fact that the $|\Lambda|\times|\Lambda|$ matrix $H_\Lambda (\theta)$ has
$|\Lambda|$  eigenvalues.
$\hfill\square$

\noindent
{\bf 2.3 The first iteration ($1^{\text {st}}$ step)}

Let 
$$
N=[N_0^C]+1, \qquad (C>1).  \tag {2.20}
$$
Let $\Lambda= [-N, N]^{d+\nu}$.
$N$ is the next scale, recall that the previous scale $N_0$ is determined by $\delta$ in (2.11).
The aim of this section is to  derive the analogue of {Lemma 2.1} for $G_{\Lambda}$.

To do that we use the estimates on $G_{\Lambda_0}$ at scale $N_0$ in {Lemma 2.1} and also
Lemma 2.2.
Let $\Lambda_0=[-N_0, N_0]^{d+\nu}+ i\subset \Lambda, i\in \Lambda$.
For a fixed $\theta$, we say that $\Lambda_0$ is {\it good} if (2.13) holds, otherwise $\Lambda_0$ is {\it
bad}.
Recall from {(2.9)}, {(2.11)} that for fixed $\theta$, at scale $N_0$, $\Lambda_0$ is bad if
$$
|n\cdot\omega+\theta+\mu_j-E|< 2 e^{-N_0^\sigma}\tag{2.21}
$$
for some $(n, j)\in \Lambda_0$, where $\mu_j$ is an eigenvalue of $H_j$.

Let $X_N$ be the set where all $\Lambda(k) =[-N_0, N_0]^d+k$, $k\in [-N, N]^d$, have the property (2.14).
Note that $X_N\subset X_{N_0}$.
So
$$
\mes \,X_N\geq 1-\frac{(2N_0+1)^{2d}}{{N'}^{2p'}} (2N+1)^d \qquad (p'>0)
\tag 2.22
$$
where 
$$
\align
N'&=[N_0^{1/\alpha}]+1, 1<\alpha<2\\
 N &=[N_0^C]+1, C> 1.
\tag  2.23
\endalign
$$
Fix $x\in X_N$. 
Assuming $2p'/\alpha- (2+C)d\gg 1$, we prove

\proclaim
{Lemma 2.3}
There exists a set 
$$
\Omega_N {\subset {(0, 1]^\nu}}, \text{\rm mes\,} \Omega_N\geq 1-e^{-N^{\frac{\sigma}{2C}}}
\tag {2.24}
$$
where $\sigma\in (0,1)$ is as in (2.11) and $C>1$ is as in (2.20), such that if $\omega\in\Omega_N$,
then for {\it any} fixed $\theta, E$, there is only {\it one} (pair-wise disjoint) bad $N_0$-box in $\Lambda =[-N,
N]^{d+\nu}$. 
Moreover $(0, 1]^\nu\backslash \Omega_N$ is contained in the union of at most $\Cal O(1)
N^{4d+\nu}$ components.
\endproclaim

\noindent
{\it Remark.}
It is crucial that  $\Omega_N$ is independent of $\theta, E$, and only depends on $x\in X_N$.

\noindent
{\bf Proof.}
Let 
$$
\matrix
&\Lambda_0= [-N_0, N_0]^{d+\nu}+ i \subset\Lambda\\
&\Lambda_0'=[-N_0, N_0]^{d+\nu} +i'\subset \Lambda 
\endmatrix
\qquad (i\not= i')\tag {2.25}
$$
be such that $\Lambda_0 \cap \Lambda_0'= \emptyset$.

Let 
$$
\align
&\Lambda_{0, j} {\text { be the projection of } } \Lambda_0 {\text  { onto } } \Bbb Z^d\\
&\Lambda_{0, n} {\text  { be the procection of } } \Lambda _0 {\text { onto } }\Bbb Z^\nu
\endalign
$$
and similarly for $\Lambda'$.

Assume both $\Lambda_0$ and $\Lambda_0'$ are bad, then there exist $(n, j)\in \Lambda, (n', j') \in \Lambda'$, such that
$$
|n\cdot\omega+\theta+\mu_j -E|< 2 e^{-N^\sigma_0}
\tag{2.26}
$$
$$
|n'\cdot\omega+\theta+\mu_{j'}-E|< 2e^{-N_0^\sigma}\tag {2.27}
$$
Subtracting {(2.27)} from {(2.26)}, we obtain
$$
|(n-n')\cdot\omega+(\mu_j-\mu_{j'})|< 4 e^{-N_0^\sigma}.\tag{2.28}
$$
Since $\Lambda_0\cap \Lambda_0'=\emptyset$.
$$
(n, j)\not= (n', j')\tag{2.29}
$$
There are 2 possibilities:

\itemitem  {$\bullet$} \ $n=n'$

In this case $\Lambda_{0, n}\cap \Lambda_{0, n'} \not= \emptyset$, so $\Lambda_{0, j}\cap\Lambda_{0, j'}=\emptyset$.
Anderson localization for $H_j$, Theorem A then implies that on $X_N$, $|\mu_j-\mu_{j'}|\geq e^{-N_0^\beta}$ for some
$$
\beta\in (0, \sigma)
\tag 2.30
$$
for all $\Lambda_{0, j}, \Lambda_{0, j'}\subset \Lambda_j, \Lambda_{0, j}\cap \Lambda_{0, j'} =\emptyset$ and any 
pair of
eigenvalues $\mu_j\in  \sigma(H_j)$, $\mu_{j'}\in \sigma(H_{j'})$ {(2.30)} 
is in contradiction with {(2.28)}.
So there can be only 1 (pairwise disjoint) bad $N_0$-box.

\itemitem {$\bullet$} \  $n\not=n'$

Let 
$$
\align
m&=n-n'\\
\lambda&=\mu_j-\mu_{j'}, \tag{2.31}
\endalign
$$
then
$$
m\in [-2N, 2N]^\nu\backslash \{0\}\tag{2.32}
$$
$\lambda$ can take on at most $(2N+1)^{2d}(2{N_0}+1)^{2d}$ different values.

So {(2.28)} corresponds to at most $\Cal O(1)N^{4d+\nu}$ inequalities in $\omega\in (0, 1]^\nu$ of the form
$$
|m\cdot\omega+\lambda|\le 4e^{-N_0^\sigma}.\tag {2.33}
$$
For each equation in (2.33), it is simple to see that the set of $\omega\in (0, 1]^\nu$ such that (2.33)
is satisfied, has one single component of measure $\leq \Cal O(1) e^{-N_0^\sigma}$.
We hence obtain the lemma for $N_1$ $N_0$ large enough.
$\hfill \square$

Let 
$$
\Lambda =[-N, N]^{d+\nu},\quad
\Cal B_x (\Lambda, E) =\{\theta|\Vert G_\Lambda(\theta, E)
\Vert\geq e^{N^\sigma}\}.\tag{2.34}
$$
From {(2.16)},
$$
\mes \, \Cal B_x(\Lambda, E)\leq e^{-\frac{N^\sigma}2} {\text { for  } }   N\gg 1.
\tag{2.35}
$$
For any $x \in X_N, X_N$ defined in (2.22, 2.14), 
using {(2.14)}, {(2.34)}  {Lemma 2.3} and resolvent expansion a
la Fr\"ohlich-Spencer, we obtain our main estimate,

\proclaim
{Lemma 2.4}
For all $\theta \in \Bbb R\backslash \Cal B_x(\Lambda_0, E)$
$$\aligned 
\Vert G_\Lambda (\omega, \theta, E)\Vert&< e^{N^\sigma},\\
|G_\Lambda (\omega, \theta, E)(m, m')|&< e^{-\gamma |m-m'|}
\endaligned\tag 2.36
$$
for all $m, m'\in \Lambda_0, |m-m'|> N/4$, $\omega\in\Omega_N\subset(0, 1]^\nu, \text{ \rm mes\,} \Omega_N\geq 
1- e^{N^{(\sigma/2C)}}$
where $\gamma'/2<\gamma<\gamma', \gamma, \sigma$ are the same as in {Lemma 2.1}, $N_0$ as defined in
{\rm (2.20)}.
\endproclaim

\noindent
{\bf 2.4 A large deviation estimate (in $\theta$) for the Green's functions at all scales}

We now build upon the estimates in Lemmas 2.1 and 2.4 to obtain estimates for Green's functions at
all scales.
In order that the bad set in $\theta$ be of small measure at larger scales, we need to let in more bad boxes
at the smaller scales.
(In {Lemma 2.1}, there is no bad box, while in {Lemma 2.4}, there is one.)
The number of bad boxes is controlled by using semi-algebraic sets as in [BGS] and {Lemma 2.3}.

Assume $0<\delta\ll 1$ is sufficiently small so that {Lemma 2.1} holds for all $N\in [N_0, N_0^{\alpha'}]$,
$\alpha'>1, N_0\gg 1$ determined by {(2.11)}.
From {Lemma 2.4}, both equations in {(2.36)} hold for all $N\in [N_0^{\alpha}, N_0^C]$ on $\Bbb
R\backslash \Cal B_{x, N}$ and $\omega \in \Omega_N$, where $C>1$ is to be determined shortly.
The probability subspace is then further restricted to be 
$$
X^{\overset {\text {def}}\to =} X_{N_0^C} =\bigcap_{N\in[N_0,
N_0^C]}X_N \tag 2.37
$$
where $X_N$  is defined similarly to (2.22).

For what is to follow, it is more convenient to slightly modify the definition and let
$$
\align
&\Cal G_x^{\gamma, \sigma} (\Lambda, E)^{\overset {\text {def}}\to =} \{\theta \in\Bbb R| \, \Vert
G_\Lambda(\theta, E)\Vert< e^{N^\sigma},\\
&\qquad|G_\Lambda (\theta, E)(m, m')|< e^{-\gamma |m-m'|}\\
&\qquad\forall m, m'\in \Lambda, |m-m'|> N/4\}\\
&\Cal B_x^{\gamma, \sigma}(\Lambda, E) \overset {\text {def}} \to = \Bbb R\backslash \Cal
G_x^{\gamma, \sigma},\tag 2.38\endalign$$
where $\gamma>0$, $0<\sigma<1, \Lambda\subset \Bbb Z^{d+\nu}$ is a cube of side length $2N+1$.

Lemma 2.1 and {Lemma 2.4} can be summarized as

\proclaim
{Proposition 2.5}
There exist $\gamma>0$, $0<\sigma<1$, such that for $0<\delta\ll 1$, $\delta_j$ satisfying (2.3), 
$0<\ve \ll 1$, there
exists $N_0$, such that for all $N\in [N_0, N_0^C](C>1)$, $\Lambda =[-N, N]^{d+\nu}\subset \Bbb Z^{d+\nu}$, cubes of
side length $2N+1$,
$$
\sup_{x \in X, E}\text{\rm mes\,} (\Cal B_x^{\gamma, \sigma}(\Lambda, E))\leq e^{-\frac{N^{\sigma}}2}\tag 2.39
$$
if $$\omega\in\Omega_x=\bigcap
_{N\in [N_0, N_0^C]}\Omega_{x, N} \subset (0, 1]^\nu$$ and $\Omega_{x, N}$ is as in Lemma
 2.3, $\mes \, \Omega_x \geq 1-e^{-N_0^{\sigma/2}}$.
\endproclaim

Let $X_{N, i}$ be defined as in (2.22) with $[-N, N]^d+ i$ replacing $[-N, N]^d$; $\Omega_{N, i}$  defined
as in {Lemma 2.3} with $\Lambda(i) =[-N, N]^{d+\nu}+i$ in place of $\Lambda=[-N, N]^{d+\nu}$. Denote by
$\text{DC}_{A,c}(M)$, the set of $\omega\in (0,1]^{\nu}$, such that (H2) is verified for $n\in[-M, M]^{\nu}$.
We now prove

\proclaim
{Lemma 2.6}
Suppose all the assumptions of {Proposition 2.5} is valid.
Let $C>10(d+\nu), 0<\sigma< 1/2, N_1= N_0^C$.
Then for all $N\in [N_1, N_1^2]$, $\Lambda$ cubes of side length $2N+1$,
let $$X=\bigcap_{N\in [N_0, N_1^2]} X_N \bigcap_{i\in [-2N_0, 2N_0]^d} X_{N_{0, i}}.$$
For any $x\in X$, let  
$$\Omega_x =\bigcap_
{N\in [N_0, N_1^2]} \Omega_{x, N} \bigcap_{i\in [-2N_0, 2N_0]^d}\Omega_{x,
N_{0, i}}
$$ 
{\rm mes } $\Omega_x\geq 1-e^{-N_0^{\sigma/2}}$.
If $\omega\in\Omega_x\cap DC_{A, c}(N^2_1)$ then
$$
\sup\limits_{x\in X, E} \text{ \rm mes }\big(\Cal B_x^{\gamma', \sigma}(\Lambda, E)\big)
\leq e^{-N^{\sigma/2}}\tag 2.40
$$
where $\gamma' =\gamma - N^{-\kappa}, \kappa =\kappa (\sigma, \gamma)>0$.
\endproclaim

\noindent
{\bf Proof.} Fix $N\in [N_1, N_1^2]$ and let
$$
\align
\Lambda &= [-N, N]^{d+\nu}, \\
T&=[-N_0, N_0]^d\times [-N, N]^\nu\tag 2.41\\
&\subset \Lambda
\endalign
$$
Let 
$$
\align
\Lambda_0&= [-N_0, N_0]^{d+\nu}\\
\Lambda_{0} (i) &=\Lambda_0+i.\tag 2.42
\endalign
$$
Define
$$
\Cal A\overset {\text {def}}\to= \bigcup\limits_{i\in [-2N_0, 2N_0]^d}
\ \Cal B_x^{\gamma, \sigma} (\Lambda_0 (i), E).\tag 2.43
$$
Since the conditions on the Green's function in (2.38) can be rewritten as polynomial inequalities by using Cramer's
rule (ratio of determinants) as in [BG, BGS], $\Cal A$ is semi-algebraic of total degree less than
$$
\align
&(2N_0+1)^{2(d+\nu)} \cdot (2N_0+1)^{2(d+\nu)}\cdot (4N_0+1)^d\\
&= \Cal O_{d, \nu}(1) N_0^{5(d+\nu)},\tag 2.44
\endalign
$$
where the first factor corresponds to the degree of each polynomial for each pair of points in a $N_0$-box, the second
is an upperbound on the $\#$ of pairs in each $N_0$-box plus the one for the Hilbert-Schmidt norm,
the third is the $\#$ of such $N_0$-boxes.
$\Cal A$ is therefore the union of at  most $\Cal O_{d, \nu}(1) N_0^{5(d+\nu)}$ intervals in $\Bbb R$ by using Theorem
1 in [Ba] (see also [BGS], where the special case we need is restated as Theorem 7.3.)

For any fixed $\theta\in \Bbb R$, let
$$
I=\{n\in [-N, N]^\nu\big| n\cdot\omega+\theta\in \Cal A\}.\tag 2.45
$$
Then for $\omega\in\Omega_N \cap DC_{A, c}$
$$
|I|\leq \Cal O_{d, \nu}(1) N_0^{5(d+\nu)}\tag 2.46
$$
by using (2.39).
This is because if there exist $n, n' \in [-N, N]^\nu$, $n\not= n'$, then for $\omega\in \text{DC}_{A,c}$
$$
|(n-n')\cdot\omega\big|\geq \frac c{N^A}\gg e^{-N ^{\sigma/2}} (c>0, A> 0).\tag 2.47
$$
Hence each interval can contain at most 1 integer $n\in [-N, N]^\nu$.

We therefore conclude that for any fixed $\theta\in \Bbb R$,
$$
\align
&\# \{i\in\Lambda\big|\Lambda_0 (i) \cap T\not= \emptyset, \Lambda_0 (i) \text { is a bad $N_0$-box}\}\\
& \leq \Cal O _{d, \nu}(1) N_0^{5(d+\nu)}\tag 2.48
\endalign
$$
where $\Lambda, T, \Lambda_0$ as defined in (2.41, 2.42).

For the $N_0$-boxes $\Lambda_0(i)(i\in\Lambda)$, such that $\Lambda_0(i)\cap T=\emptyset$, we use Lemma 2.3 and (2.3)
to conclude that for $\omega\in\Omega_N, (N\in [N_1, N_1^2])$, $\exists \, i_0 \in\Lambda\backslash \Bbb T$, such that
$\forall i\in\Lambda$, such that $\Lambda_0(i) \cap T=\emptyset$, if $\Lambda_0(i) \cap \Lambda_0(i_0)=\emptyset$,
then
$$
\Lambda_0(i) \text { is a good $N_0$-box}.\tag 2.49
$$
We now introduce an intermediate scale $\bar N$:
$$
\log N_0< \log\bar N< \log N.\tag 2.50
$$
Let $\bar C\in \big(10(d+\nu), C\big), \bar N=[N_0^{\bar C}]$
$$
\bar\Lambda= [-\bar N, \bar N]^{d+\nu}.\tag 2.51
$$
Let $\Cal I$ be the set defined in (2.48).

We say $\bar \Lambda (i) (i\in \Lambda)$ is {\it good} if
$$
\bar\Lambda(i) \cap (\Cal I\cup \Lambda_0(i_0))=\emptyset, \tag 2.52
$$
where $\Lambda_0 (i_0)$ is as in (2.49), otherwise it is {\it bad}.

Let $\Cal F$ be a family of pairwise disjoint bad $\bar N$-boxes in $\Lambda =[-N, N]^{d+\nu}$. (2.48, 2.49) imply that
$$
\align
\#\Cal F&\leq \Cal O_{d, \nu}(1) N_0^{5(d+\nu)} +2^{d+\nu}\\
&= \Cal O_{d, \nu}(1) N_0^{5(d+\nu)}.\tag 2.53
\endalign
$$
If $\bar \Lambda(i)$ is good, then
$$
\forall j\in\bar\Lambda(i), \exists\,j'\in \bar\Lambda (i) \text { such that $\Lambda_0(j')\subset
 \bar \Lambda(i)$}
$$
and 
$$
\dist \big(j, \partial_*\Lambda_0(j')\big) \geq N_0,\tag 2.54
$$
where $\partial_*\Lambda_0(j')$ is the interior boundary of $\Lambda_0(j')$ relative to $\bar\Lambda(i)$:
$$
\partial_*\Lambda_0 (j')=\{ z\in \Lambda_0(j')|\exists \, z' \in\bar\Lambda(i)\backslash \Lambda_0(j'), |z'- z|=1\}.\tag
2.55
$$
An easy resolvent expansion (see e.g. Lemma 2.2 and proof of Corollary 4.5 in [BGS]) then shows that
$$
|G_{\bar \Lambda(i)}(m, m')|< e^{-\gamma|m-m'|+CN_0}\tag 2.56
$$
for all $m, m' \in\bar \Lambda(i), |m-m'|> N_0/4$.

For $N\in [N_1, N_1^2]=[N_0^C, N_0^{2C}]$,
$$
\# \Cal F< N^\sigma\tag 2.57
$$
with
$$
\frac{5(d+\nu)}{2C}< \sigma < \frac{5(d+\nu)}C\tag 2.58
$$
for 
$$
C>10(d+\nu), \sigma< 1/2.\tag 2.59
$$

For all $C> 10(d+\nu)$, we can choose $\bar C$ satisfying (2.51), so that for all $N\in [N_1, N_1^2]$, we obtain
(2.40) from the estimates at scale $\bar N$ by applying Lemma 2.1 of [BGS] with a single step iteration.
(This is possible because $\sigma<1/2$, so $\exists \, \alpha>2$ such that $\alpha \sigma< 1$, where $\alpha$ is the
geometric expansion factor.
See the first inequality of (2.12) of [BGS].)

The measure estimate in $\theta$ is supplied by Lemma 2.2.

We do not repeat the details of this iteration, except noting the following small variations:

\itemitem {$\bullet$}  Because of the apriori estimate in Lemma 2.2, which holds at all scales, we only need to
estimate Green's functions for cubes in $\Bbb Z^{d+\nu}$.

\itemitem {$\bullet$}  To eliminate $G_\Lambda(i, i'), i, i'\in \Lambda, |i-i'|> N/4$, we make an exhaustion
$\{S_j(i)\}^\ell_{i=0}$ of $\Lambda$ of width $2\bar N$ centered at $i$ (as in [BGS]):
$$
\align
S_{-1} (i) &\overset \text{def}\to = \emptyset\\
S_0(i) &\overset \text{def}\to =\bar\Lambda (i) \cap\Lambda\\
S_j(i)&\overset\text{def}\to = \bigcup_{k\in \delta_{j-1}(i)} \bar\Lambda(k) \cap \Lambda
\text { for } 1\leq j\leq \ell\tag 2.60
\endalign
$$
where $\ell$ is maximal such that $S_\ell (i)\not= \Lambda$.

\itemitem{$\bullet$} 
In the iteration, we need to estimate
$ G_A(m', m)$,  where $A=S_j\backslash S_{j'}$, for some $j'<j\leq\ell$, is an annulus.
$m'\in  \partial_{**}S_{j'}, \partial_{**}S_{j'}$, is the exterior boundary of $S_{j'}$, relative to $\Lambda$:
$$
\partial_{**}S_{j'}=\{ z| z\in \Lambda\backslash S_{j'}, \exists\, z' \in S_{j'}, \quad |z-z'|=1\}\tag 2.61
$$
is concave;
$m\in\partial_{*} S_j$, $\partial_{*}S_j$ is the interior boundary of $S_j$ as defined in (2.55)
$\partial_*S_j$ is convex.

\itemitem{$\bullet$} \
When $A$ is {\it good} (for the precise definition, see Lemma 2.2 of [BGS]).
We estimate $G_A(m', m)$ using $\bar \Lambda$ cubes, which are all ``good''.
We always start the resolvent expansion from $m\in \partial_*S_j$, which  for the property that $\forall m\in
\partial_*S_j$, $\exists \, m'' \in A$, such that $\bar\Lambda(m'') \in A$ and $\dist \big(m, \partial_* \bar\Lambda
(m'')\big)\geq N_0$.
For the last term in the expansion, we use the apriori estimate in Lemma 2.2.
This way we avoid having to estimate Green's functions in regions of  the form $\bar\Lambda(z)\cap\Lambda$.
We obtain exponential decay as in (2.56).

\itemitem {$\bullet$} 
When $A$ is {\it bad}, we resort to Lemma 2.2.
From (2.60), we need (2.16) to hold for at most
$$
\Cal O_{d, \nu}(1)\frac N{\bar N}\cdot N^{d+\nu}\leq\Cal O_{d, \nu} (1) N^{d+\nu+1}\tag 2.62
$$
number of annuli.
Combining (2.16) with (2.62) we obtain the estimate in measure in (2.40).
$\hfill\square$

Using Proposition 2.5 and Lemma 2.6, we obtain the main estimate of this section by induction:

\proclaim 
{Proposition 2.7}
There exist $\gamma>0, 0<\sigma< 1/2, N_0\in \Bbb N, X\subset\Bbb R^{\Bbb Z^d}$,
$$
\align
&\text{\rm mes\,} X\geq 1 -N_0^{-1}\tag 2.63\\
&\Omega_x\subset (0, 1]^\nu,
\text{\rm mes}\, \Omega_x\geq 1-e^{ - N_0^{\sigma/2}},\tag 2.64\\
&\text {such that for $x\in X$, $\omega\in\Omega_x\cap DC_{A, c}$},\\
&\text{$0<\delta\ll 1$, $\delta_j$ satisfying (2.3), for all
$N\geq N_0$},\\
&\text{$\Lambda =[-N, N]^{d+\nu}+i$, all $i\in [-2N, 2N]^d$},\\
&\sup_{x\in X, E}\text{ \rm mes\,} (\Cal B_x^{\gamma, \sigma}(\Lambda, E))\leq e^{-N^{\sigma/2}}.\tag 2.65
\endalign
$$
\endproclaim

\noindent
{\bf Proof.}
In view of the proof of Lemma 2.1, in particular (2.22) and the proof of Lemma 2.6:
$$
\mes(\Bbb R^{\Bbb Z^d}\backslash X)\leq \sum^\infty_{L_i=N_0} \ \sum_{L_i^c\leq L_{i+1}\leq L_i^{2C}}
L_{i-1}^{-2p'}(2L_i+1)^{2d}(6L_{i+1}+1)^d\tag 2.66
$$
$$
\align
&L_{i-1}= L_i/\alpha \qquad (1<\alpha<2)\\
&C> 10(d+\nu)\tag 2.67
\endalign
$$
where the factor 6 in the last factor of (2.66) comes from the fact that at each scale $L$, we need estimates for all
cubes $[-L, L]^d+i, i\in[-2L, 2L]^d$, see (2.43).

(2.63) is satisfied if
$$
\bar N_0^{2p'/\alpha} N_0^{2d} N_0^{2C} N_0^{2Cd}\leq N_0^{-2},\tag 2.68
$$
which leads to
$$
p'>\alpha (C+1) (d+1).
$$
In view of (2.67)
$$
p'> 42(d+\nu)d.\tag 2.69
$$
(2.64) is verified by removing a set $\Omega_L$ at each scale
$L$ similar to Lemma 2.3.
Because of the decay property of $\delta_j$ in (2.3), the RHS of (2.26, 2.27) are replaced by $\Cal O(1) e^{-L'}$
($\Cal O(1)$ is the same for all scales), where $L'$ is the previous scale, see proof of Lemma 2.1.
Summing over the scales, we obtain (2.64).
$\hfill\square$

\noindent
\heading
{3. Exponential Decay of the Green's function of Schr\"odinger operator at fixed $E$ and $\theta$}
\endheading

This section is in some sense a minor image of the previous section.
Here we study the operator $H(\theta)$ defined in (2.5) for fixed $\theta$, but we are allowed to ``move'' the random
variables $x\in \Bbb R^{\Bbb Z^d}$.
Without loss of generality, we set $\theta=0$ and study
$$
H=\delta_j\Delta_n+n\cdot\omega+\epsilon\Delta_j+V_j\tag 3.1
$$
on $\ell^2(\Bbb Z^{d+\nu})$.
We first prove the analogue of Lemma 2.2.

\proclaim
{Lemma 3.1}
Let $E\in I$, an interval of length $\Cal O(1)$.
Let $\Lambda\subset\Bbb Z^{d+\nu}$ be a finite set.
Then
$$
\text{\rm mes\,} \{x \big|\text {\rm dist\,}
\big(E, H_\Lambda (x)\big) \leq \kappa\} \leq C\kappa |\Lambda| \ \Vert\tilde g\Vert_\infty.\tag 3.2
$$
\endproclaim

\noindent
{\bf Proof.}
Let $\Cal N(E, x)$ be the $\#$ of eigenvalues of $H_\Lambda\leq \lambda$
$$
\align
&\mes \{x\big|\dist \big(E, H_\Lambda(x)\big) \leq \kappa\}\\
&\leq \int\big(\Cal N(E+\kappa, x) -\Cal N(E-\kappa, x)\big) \prod_{i\in\Lambda_j} \tilde g(v_i) dv_i\\
&=\iint^{E+\kappa}_{ E-\kappa} \frac d{d\lambda}\Cal N(\lambda, x) d\lambda \prod_{i\in\Lambda_j} \tilde g
(v_i) dv_i\\
&=\int^{E+\kappa}_{E-\kappa} \sum_{i\in\Lambda_j}\frac\partial{\partial v_i} \Cal N(\Lambda, x)\prod_{i\in
\Lambda_j} \tilde g(v_i)dv_i\\
&\leq C\kappa |\Lambda|\Vert \tilde g\Vert_\infty
\tag 3.3
\endalign
$$
where $\Lambda_j =\Lambda \cap \Bbb Z^d$ and each $\partial/{\partial v_i}$ is seen as a rank $|\Lambda\cap (\Bbb
Z^\nu+i)|$ perturbation as usual.
~$\hfill\square$

We define the good and bad sets in analogy with (2.38):
$$
\align
\Cal G_\theta^{\gamma, \sigma}(\Lambda, E) &\overset{\text{def}}\to = 
\{x\in\Bbb R^{\Bbb Z^d}\big| \ \Vert G_\lambda(x, E) \Vert < e^{N^\sigma},
\\
&|G_\Lambda (x, E) (m , m')|< e^{-\gamma |m-m|} \\
& \forall m, m' \in\Lambda, |m-m'| >N/4\}
\endalign
$$
$$
\Cal B_\theta^{\gamma, \sigma} (\Lambda, E)\overset\text{def}\to = \Bbb R^{\Bbb Z^d}\backslash \Cal G_\theta^{\gamma,
\sigma} (\Lambda, E)\tag 3.4
$$
where $\gamma>0, 0< \sigma< 1, \Lambda =[-N, N]^{d+\nu} +i$ for some $i\in\Bbb Z^{d+\nu}$.

\proclaim
{Lemma 3.2} There exist $\gamma>0, 0<\sigma<1, 0<\epsilon_0 \ll 1, N_0(\sigma), p>0$,
such that for all $0<\epsilon <\epsilon_0, N\geq N_0$, there exists $0<\delta(N)\ll 1$, such that
$$
\sup_{\theta, E} \text{\rm mes\,} \big(\Cal B_\theta^{\gamma, \sigma}(\Lambda, E)\big)
\leq N^{-p}.\tag 3.5
$$
\endproclaim

\noindent
{\bf Proof.}
We prove (3.5) by perturbation.
Let
$$
H_{\Lambda, 0} \overset{\text{def}}\to= n\cdot \omega+\epsilon\Delta_j+V_j\tag 3.6
$$
on $\ell^2(\Lambda)$ with Dirichlet boundary conditions as defined in (2.4),
$$
H_j\overset{\text{def}}\to = \epsilon\Delta_j+V_j\tag 3.7
$$
on $\ell^2(\Lambda\cap \Bbb Z^d)$.

For $\epsilon$ sufficiently small, $H_j$ has Anderson localization, (cf. appx.)
So for any fixed $\lambda$ starting at some scale $N'$, with log $N'<\log N$, using Theorem A and Lemma 3.1,
$$
|(H_j-\lambda)^{-1} (j, j')|< e^{-\gamma|j-j'|}\tag 3.8
$$
$$
\text{ for all $j, j'\in \Lambda\cap \Bbb Z^d, |j-j'|> N/4$}
$$
with probability $\geq 1-\Cal O_d(1) {N^{2d}}/{N'{}^{2p'}}-e^{-N^\sigma}$,
$0<\sigma< 1, p'> 0$ to be specified, $N'=N^{1/\alpha}(1<\alpha< 2)$ and where we used Lemma 3.1.

Let 
$$
\lambda =E+n\cdot\omega,\quad n\in\Lambda\cap \Bbb Z^\nu\tag 3.9
$$
and summing over the probability, we obtain that
$$
|(H_{\Lambda, 0} -E)^{-1} (j, j')|< e^{-\gamma|j-j'|}\tag 3.10
$$
for all $j, j'\in \Lambda \cap \Bbb Z^d, |j-j'|> N/4$ with probability
$$
\geq 1-\Cal O_d(1) \frac{N^{2d}}{N'{}^{2 p'}} -\Cal O_\nu(1) N^\nu \cdot e^{-N^\sigma}.
$$

We remark here that the set defined by (3.10), $X_{N'}\subset X_{N_0}$, $X_{N_0}$ as in (2.14), $\mes X_{N'}\leq\mes
X_{N_0}$, because of the need for Lemma 3.1, as we do not have $\theta $ at our disposal.

To obtain exponential decay of $G_\Lambda$, we use the resolvent expansion:
$$
\align
(H_\Lambda -E)^{-1} (n, j, n', j') &= (H_{\Lambda, 0} -E)^{-1} (n, j; n', j')\\
&+\bigg(\sum_{k=1}^K[(H_{\Lambda, 0}-E)^{-1} \delta_j\Delta_n]^k (H_{\Lambda, 0}-E)^{-1}\bigg) (n, j; n', j')\\
&+ \big([(H_{\Lambda, 0}-E)^{-1}\delta_j\Delta_n]^{K+1} (H_\Lambda -E)^{-1}\big) (n, j, n', j')\tag 3.11
\endalign
$$
where $K=\Cal O(n-n')$.
Using (3.10, 3.2, 2.3), we arrive at (3.5) for
$$
\delta=\delta(N) \ll 1, N'=[N^{1/\alpha}], 1< \alpha<2, p=\frac{2p'}\alpha -2d\tag 3.12
$$
(assuming $p'$ large enough),
$N$ large enough and a $\gamma$ which is slightly smaller than that in (3.10).
$\hfill\square$

Lemma 3.2 is the analogue of Proposition 2.5.
We now prove
\proclaim
{Lemma 3.3}
Suppose $\delta$ is such that Lemma 3.2 holds for all $N\in [N_0, N_0^C]$ with $C> 10(d+\nu)$.
Let $N_1 =N_0^C$.
Then for all $N\in [N_1, N_1^2]$,
$ \Lambda =[-N, N]^{d+\nu}+i, i\in\Bbb Z^{d+\nu}$, if $\omega\in DC_{A, c}(2N)
$
then 
$$
\sup_{E, \theta}\text{\rm mes\,} \big(\Cal B_\theta^{\gamma, \sigma}(\Lambda, E)\big) \leq N^{-p}\tag 3.13
$$
where $\gamma'=\gamma -N^{-\kappa}, \kappa =\kappa(\sigma, \gamma)> 0$, $p$ is the same as in Lemma 3.2, (3.12).
\endproclaim

\noindent
{\bf Proof.}
This is almost a mirror image of the proof of Lemma 2.6.
Hence we will only point out the small differences.

Fix $N\in [N_1, N_1^2]$ and let
$$
\align
\Lambda& = [-N, N]^{d+\nu}\\
T(i) &= \{[-N_0, N_0]^2 +i\}\times [-N, N]^\nu \qquad (i\in \Lambda\cap \Bbb Z^d)\\
\Lambda_0&=[-N_0, N_0]^{d+\nu}\\
\Lambda_0(i)&= \Lambda_0+i \qquad (i\in \Lambda).\tag 3.14
\endalign
$$

\itemitem {$\bullet$} \ For a given $T^{(i)}, i\in \Lambda\cap \Bbb Z^d$ 
$$
\text {Prob\,}\{\Lambda_0(i+j) \text{ are good, } \forall j\in [-N, N]^\nu\}\geq 1-\Cal O_\nu(1)
\frac{N^\nu}{N^{p}_0}\tag 3.15
$$
by using Lemma 3.2.

\itemitem {$\bullet$} For any $x\in X_{N_0}' \subset X_{N_0}$, for a given $T(i)$, from semi-algebraic considerations
as in (2.43-2.48) for $\omega\in DC_{A, c} (2N)$
$$
\align
&\# \{i'\in [-N, N]^{d+\nu}| \Lambda_0(i')\cap T(i)\not= \emptyset, \Lambda_0(i') \text { is a bad $N_0$-box}\}\\
&\leq \Cal O_{d, \nu} (1) N_0^{5(d+\nu)}\tag 3.16
\endalign
$$
as in (2.48).

For a given
$$
T(i) =\{[-N_0, N_0]^d+i\} \times [-N, N]^\nu,
$$
$i\in\Lambda\cap\Bbb Z^d$, we say $T(i)$ is {\it bad} if the condition in (3.15) is violated, i.e., if $\exists\,
 j\in [-N, N]^\nu$, such that $\Lambda_0(i+j)$ is a bad $N_0$-box.
Let $J$ be an even integer.
From (3.15)
$$
\align
&\text {Prob  \{there are at least $J$ pairwise disjoint bad $T(i)$ in $\Lambda\}$}\\
&\leq [\bigg(\Cal O_\nu(1) \frac{N^\nu}{N^p_0}\bigg)^2 \cdot \Cal O_d(1)N^{2d}]^{J/2}\\
&= \Cal O_{d,\nu}(1) \bigg(\frac{N^{d+\nu}}{N^p_0}\bigg)^J\\
&= \Cal O_{d,\nu}(1) N^{-(\frac p{2C} -d-\nu)J}\tag 3.17
\endalign
$$
assuming 
$$
p>2C(d+\nu)> 20(d+\nu)^2.\tag 3.18
$$

For a given $C$, let $p$ be such that $ p/2C-d-\nu=1$ and choose $J= p+1$ or $J=p+2$, so that $J$ is even.
So
$$
(3.17) \leq \Cal O_{d, \nu} (1)N^{-(p+1)}.\tag 3.19
$$

\itemitem {$\bullet$}  Using (3.16, 3.19, 3.2), introducing another scale $\bar N=N_0^{\bar C}$, $\bar C\in
(10(d+\nu), C)$, we obtain Lemma 3.2 just like the way we obtained Lemma 2.6.

Iterating using Lemmas 3.2, 3.3, we arrrive at the main estimate of this section.

\proclaim
{Proposition 3.4} There exist $\gamma>0, 0< \sigma <1/2$, $p>0$, such that for $\omega\in D C_{A, c}, 0<\epsilon \ll 1,
0<\delta \ll 1$, $\delta_j$ satisfying (2.3), there exists $N_0(\epsilon, \delta)$, such that for all $N\geq N_0,
\Lambda\subset \Bbb Z^{d+\nu}$, cubes of side length $2N+1$
$$
\sup_{\theta, E} \text{\rm mes\,} \big(\Cal B_\theta^{\gamma, \sigma}(\Lambda, E)\big) \leq N^{-p}.\tag 3.20
$$
\endproclaim
\bigskip

\noindent
\heading
{4. The elimination of $E$ and frequency estimates}
\endheading

The goal of this section is to transform the large deviation estimates in $\theta$ in (2.63) into estimates in
$\omega$ and to eliminate the dependence of the singular set on $E$ in the process.
(Recall that $\theta$ is an auxilary variable that was not in the original problem (2.2 ). )
This is needed to prove Anderson localization.
We need two scales $N,\bar N$, $\log \log\bar N\ll \log N$.
Let $\Lambda_N(i)=[-N,N]^{d+\nu}+i, i\in \Bbb Z^d, \Lambda_{\bar N}= [-\bar N, \bar N]^{d+\nu}, \Omega_N(i)$ the set
of $\omega$ as defined in Lemma 2.3 for the box $\Lambda_N^{(i)}$ and $DC_{A, c}(N)$, the set of $\omega$ satisfying
(H2) for $n\in [-N, N]^\nu$.

For a fixed $x\in X', X' \subset X$ as in (2.63), define $\Cal B_x(N, \bar N)\subset (0, 1]^\nu\times \Bbb R$ as
$$
\align
\Cal B_x(N, \bar N)&=\bigcup_{\{i\in\Bbb Z^d|\Lambda_N(i) \cap \Lambda_{\bar N}\not= \emptyset\}} 
\{(\omega, \theta)\in (0, 1]^\nu \times\Bbb R|\exists \, E, \text { such that }\\
&\Vert\big(H_{\bar N}(\omega, 0) -E\big) ^{-1} \Vert \geq e^{\bar CN} \text { and } \theta \in\Cal B_x^{\gamma,
\sigma} \big(\Lambda_N(i), E\big)\}\tag 4.1
\endalign
$$
and $S_x(N) \subset (0, 1]^\nu\times \Bbb R$ as
$$
\align
S_x(N)=\bigcup_{\bar N\asymp  N^C} &\big\{\Cal B_x(N, \bar N)\bigcap (\cap_{\{i\in\Bbb Z^d|\Lambda_N(i) \cap
\Lambda_{\bar N}\not=\emptyset\}} (\Omega_N(i) \times\Bbb R))\big\}\\
& \bigcap (DC_{A, c}(\bar N+N)\times\Bbb R)\tag 4.2
\endalign
$$
where $\Omega_N(i)$ is the set of $\omega$ defined as in Lemma 2.3 for the box $[-N, N]^{d+\nu}+i$.

In view of (4.1, 4.2), at each scale $N$, we need estimates for all cubes $[-N, N]^d+i, i\in [-\bar N, \bar N]^d$,
where $\bar N\asymp N^C$.
Similar considerations as in (2.67-2.69) lead to $p'>422 d(d+\nu)^2$
$$
\mes X' \geq 1-N_0^{-1}\tag 4.3
$$
$N_0$ as in Proposition 2.7.

Likewise $\Omega_x$ of Proposition 2.7 is reduced to $\Omega_x' \subset\Omega_x$ accordingly,
$$
\mes\Omega_x' \geq 1-e^{-N_0^{\sigma/2}}\tag 4.4
$$
if $N_0\gg 1$.
For simplicity, we now drop the prime:
$$
X{\overset{\text {def}}\to =} X', \Omega_x \overset\text{def}\to = \Omega_x'.
$$

 \proclaim
{Lemma 4.1} Let $N, \bar N\in \Bbb N$ be such that $\bar N\asymp N^C (C>1)$.
For any $x\in X$,
$$
\mes S_x(N) < e^{-\frac 12 N^{\sigma/2}} \qquad 0<\sigma < 1/2;\tag 4.5
$$
Moreover for any $\theta\in\Bbb R$, the section
$$
S_{x,\theta} (N) =\{\omega\in (0, 1]^\nu |(\omega, \theta)\in S_x(N)\}
$$
is a union of at most $N^{7C(d+\nu)\nu}$ components.
\endproclaim

\noindent
{\bf Proof.}
For a given $\bar N$, fix
$$
\omega \in\bigcap_{\{ i\in\Bbb Z^d|\Lambda_N(i) \cap \Lambda_{\bar N}\not = \emptyset\}}
\Omega_N(i)\bigcap DC_{A, c}(\bar N+N).\tag 4.6
$$
Let
$$
\lambda \in \sigma \big(H_{\bar N} (\omega, 0)\big).$$
Then
$$
\align
&\mes \{ \bigcup_{\{i\in\Bbb Z^d|\Lambda_N(i)\cap \Lambda_{\bar N}\not=\emptyset\}}
 \ \bigcup_{ \lambda\in \sigma(H_{\bar N}(\omega, 0))} \Cal B_x^{\gamma, \sigma}(\Lambda_N(i), \lambda)\}\\
& \leq \Cal O_{d, \nu} (1) \bar N^d\cdot \bar N^{d+\nu}\cdot e^{-N^{\sigma/2}}\\
&\leq\Cal O_{d, \nu}(1) e^{-\frac 12 N^{\sigma/2}}\tag 4.7
\endalign
$$
for $\bar N\asymp N^C$, where we used Proposition 2.7.

Let $E$ be such that 
$$
\Vert (H_{\bar N} (\omega, 0)-E)^{-1}\Vert \geq e^{\bar C N},
$$
then $\exists \, \lambda \in \sigma(H_{\bar N}(\omega, 0))$ such that
$$
|E-\lambda|\leq e^{-\bar CN}.\tag 4.8
$$
Using the resolvent equation, we have
$$
(H_{\Lambda_N(i)} -\lambda)^{-1}=(H_{\Lambda_N(i)} -E)^{-1}
+(E-\lambda)(H_{\Lambda_N(i)} -\lambda)^{-1} (H_{\Lambda_N(i)}-E)^{-1}.\tag 4.9
$$
(4.9) and Lemma 2.2 then imply that
$$
\Cal G^{\gamma, \sigma}_{x}(E) \supset \Cal G_x^{\gamma', \sigma'} (\lambda)\tag 4.10
$$
with  $0 <\gamma'<\gamma, 0<\sigma'< \sigma$.
Using (4.10) and taking the union over $\bar N\asymp N^C$, we obtain (4.5) with $\gamma, \sigma$ slightly smaller
than that in (4.7).

To prove the second statement, we need to bound the degree of $S_{x,\theta} (N)$ for a fixed $\theta$.
The conditions in (4.2) can be reexpressed in polynomial inequalities by using Hilbert Schmidt norm and Cramer's
rule as before.
The sets $\Omega_N(i)$ are defined by polynomial (monomial) inequalities as in (2.28, 2.33), similarly for
$DC_{A,c}(\bar N+N)$. So $S_x(N)$ is semi-algebraic.

Using a special case of Theroem 1 in [Ba] as stated in Theorem 7.3 of [BGS], we obtain that for any fixed $\theta$,
$S_{x, \theta}(N)$ is the union of at most
$$
\align
&\Cal O_{d, \nu}(1) \{\bar N^{2(d+\nu)} (N^{2(d+\nu)}+N^{2(d+\nu)}
+N^{2(d+\nu)})\cdot \bar N^{d}\cdot N^c\cdot \bar N^{d+\nu}\}^\nu\\
& \leq \Cal O_{d, \nu}(1) N^{7C(d+\nu)\nu}
\endalign
$$
connected components.
$\hfill\square$

We need one more lemma, before transfering the estimate in $(\omega, \theta)$ in (4.5) into an estimate in $\omega$
only.

\proclaim
{Lemma 4.2}
Let
$$
S(N)\subset \big((0, 1]^\nu \cap DC_{A, c}(N)\big) \times \Bbb R
$$
be a set with the properties:

\itemitem {$\bullet$} For each $\theta\in\Bbb R$, the section
$$
S_\theta=\{\omega\in \{(0, 1]^\nu \cap DC_{A, c}(N) \}|(\omega, \theta)\in S\}
$$
is a union of at most $M$ components

\itemitem {$\bullet$} $\text{\rm supp\,} S\subset (0, 1]^\nu \times[-N_0, N_0]
$.

Let $K\gg MN_0$.
Then
$$
\align
&\text{\rm mes\,}\{\omega\in \{(0, 1]^\nu \cap DC_{A, c}(2K)\}\big|(\omega, \ell\cdot \omega)\in S(N)
\text{ for some $\ell, |\ell|\sim K$}\}
\\
&\leq \Cal O(1)(K^\nu MN_0 \big(\text{mes}\, S(N)\big)^{1/3} +MN_0K^{-1}).\tag 4.11
\endalign
$$
\endproclaim

\noindent
{\bf Proof.}
We use a similar strategy as in the proof of Lemma 6.1 in \cite{BG}.
$$
\align
&\mes \{\omega\in \{ (0, 1]^\nu \cap DC_{A, c}(2K)\}\big|(\omega, \ell\cdot\omega)\in S(N) \text { for some }
\ell, |\ell| \sim K\}\\
&\leq \sum_{\ell, |\ell|\sim K} \int_{(0, 1]^\nu} \chi_{S(N)} (\omega, \ell\cdot\omega) d\omega.\tag 4.12
\endalign
$$
Let
$$
\ell\cdot\omega\overset\text {def}\to = \theta =|\ell|\omega_\ell, \tag 4.13
$$
where $\omega_\ell$ is the projection of $\omega$ in the $\ell$ direction.
Let $\omega_\ell^\bot$ be the orthogonal component.
So
$$
(4.12) \leq\sum_{\ell, |\ell|\sim K} \frac 1{|\ell|} \int_{-N_0}^{N_0} d\theta \int d\omega_\ell^\bot \chi_{S(N)}
\bigg(\bigg(\frac\theta{|\ell|}, \omega^\bot_\ell\bigg), \theta\bigg).\tag 4.14
$$
Fix $\theta$ and bound
$$
\# \{\ell| \ |\ell| \lesssim K \text { and } \omega=\bigg(\frac\theta{|\ell|}, \omega^\bot_\ell\bigg) \in S_\theta (N)\}.
$$
We distinguish two cases:

\itemitem {$\bullet$} $|S_\theta|>\gamma$ \hfill (4.15)

\itemitem{$\bullet$} $|S_\theta|\leq \gamma$\hfill (4.16)

where $0<\gamma \ll 1$ is to be specified.

\itemitem {$\bullet$}  If $|S_\theta|>\gamma$, then the contribution to (4.14) is bounded by
$$
\align
&K^{\nu-1} \{\mes  \theta\in[-N_0, N_0]| \ |S_\theta|>\gamma\}\\
&\leq \Cal O_\nu (1) K^{\nu-1} \gamma^{-1}|S|.\tag 4.17
\endalign
$$

\itemitem {$\bullet$} \ Assume (4.16) and $|\theta|< aK^{-
A}$ ($a>0$ to be specified), the contribution to (4.14) is bounded by
$$
\Cal O_\nu (1) K^{\nu-A-1} a.\tag 4.18
$$

\itemitem{$\bullet$} \ Assume (4.16), $|\theta|\geq aK^{-A}$ and moreover
$$
\exists \, \ell, \ell', \ell\not=\ell', \text { such that $\exists\,  \omega, \omega'$ in the same component
of $S_\theta$}\tag 4.19
$$
Since $\theta$ is fixed
$$
|\ell\cdot \omega -\ell'\cdot\omega'|=0\tag 4.20
$$
$$
\ell\cdot\omega -\ell'\cdot\omega' = (\ell -\ell') \cdot\omega+\ell'\cdot (\omega-\omega')
$$
$$
|(\ell -\ell')\cdot\omega|\geq \frac a{(2K)^A}
$$
for $\omega \in DC_{A, a}(2K)$.
So (4.19) implies
$$
|\ell'\cdot (\omega-\omega')|\geq \frac a{(2K)}.
$$
Hence
$$
|\omega_\ell -\omega_{\ell'}|\geq \frac a{(2K)^{A+1}}\tag 4.21
$$
(4.16, 4.21) imply that the contribution of (4.19) to (4.14) is bounded by
$$
\align
&\Cal O(1) K^{\nu-1} N_0M\gamma K^{A+1} a^{-1}\\
&=\Cal O(1) K^{\nu+A}N_0 Ma^{-1}\gamma.\tag 4.22
\endalign
$$

The contribution from the negation of (4.19) is bounded by
$$
\Cal O(1)K^{-1} N_0M.\tag 4.23
$$
Summing over (4.17, 4.18, 4.22, 4.23) and taking $\gamma=|S|^{2/3}, a=|S|^{1/3} K^A$, we obtain the lemma.
$\hfill\square$

Combining Lemmas 4.1 and 4.2, we arrive at the conclusion of this section.
Let $\bar N\asymp N^C$,
$$
\tilde\Omega_{N, \bar N}=\bigcap_{\{i\in\Bbb Z^d|\Lambda_N(i) \cap \Lambda_{\bar N}\not=\emptyset\}}
\Omega_N(i)\tag 4.24
$$
$\Omega_N(i)$ as defined in Lemma 2.3 for the box $[-N, N]^{d+\nu}+i$ we have

\proclaim
{Lemma 4.3}
Let $N\in \Bbb N$ be sufficiently large.
Fix $x\in X$.
Let
$$
\tilde \Omega_N= \tilde\Omega_{N, \bar N} \cap DC_{A, c} (\bar N+N)
$$
be the set such that

\noindent
(4.25) There is $\bar N\asymp N^C, \ell\in \Bbb Z^\nu, |\ell|\sim N^\tau (\tau>7C(d+\nu)\nu)$
and $E$ such that
$$
\Vert\big(H_{\bar N}(\omega, 0)-E\big)^{-1} \Vert\geq e^{\bar CN}\qquad (\bar C>0)\tag 4.26
$$
and there is $i\in \Bbb Z^d$, $\Lambda_N(i) \cap \Lambda_{\bar N}\not=\emptyset$.
Such that
$$
|(H_{\Lambda_N}(i)-E)^{-1}(\omega, \ell\cdot\omega)(m, m')|>e^{-\gamma|m-m'|}\tag 4.27
$$
for some $m, m' \in\Lambda_N(i), |m-m'|> N/4$, 
$$
\text{\rm mes\,}\tilde\Omega_N\leq N^{-\tilde q}\qquad 0< \tilde q <(\tau -7C(d+\nu)\nu-1).\tag 4.28
$$
\endproclaim

\noindent
{\it Remarks.}
\itemitem{$\bullet$} \ $\ell$ could be taken larger, e.g. $|\ell|\sim N^{\log N}$ as in \cite{BG, BGS}.
But in view of the probability estimate for the random part, which is only polynomial coming from \cite
{vDK}.
We take $|\ell| \sim N^{\tau}$, assuming $p\gg 1$ ($p$ as in (3.12)).

\itemitem {$\bullet$} \ The probability estimate for the random part can be improved to subexponential by allowing
more bad boxes. But for now, we leave it as it is.
\bigskip
\heading
5. Proof of Anderson Localization for the Schr\"odinger operator
\endheading

We now prove Anderson localization, i.e., pure point spectrum with exponentially decaying eigenfunctions, for $H$
defined in (2.2).

\proclaim
{Theorem}
There exists $\epsilon_0$, such that $\forall\, 0<\epsilon <\epsilon_0, \eta>0,$ given a bounded interval $I\subset\Bbb
R$, $ \exists\, \delta_0(\eta, \epsilon)$, such that $\forall 0<\delta<\delta_0(\eta, \epsilon)$
$$
\exists\, \tilde X_{\eta,\epsilon,W} \subset\Bbb R^{\Bbb Z^d},
\,\text{\rm mes}\tilde X_{\eta, \epsilon,W}\geq 1-\eta,$$ such that

$$\forall x\in\tilde X_{\eta, \epsilon, W}, \exists\, \Omega_{x, \eta, \epsilon, W} \subset (0, 1]^\nu,
$$
$\text{\rm mes\,}\Omega_{x, \eta, \epsilon, W}\geq 1-\eta$, such that
$$
\forall \omega\in \Omega_{x, \eta,\epsilon, W},  \delta_j \text { satisfying (2.3)}
$$
$H$ has Anderson localization in $I$.
\endproclaim

We need the analogue of Lemma 4.3.

\proclaim
{Lemma 5.1} Let $N\in \Bbb N$ be sufficiently large.
Fix $\theta =0$.
Let $\omega\in DC_{A, c} (N)$.
$$
 \text {There is  } \bar N\asymp N^C, \ell\in \Bbb Z^{d+\nu}, |\ell|\sim N^\tau,\, E\tag 5.1
$$
such that
$$
\Vert\big(H_{\bar N}(\omega, 0) -E\big)^{-1}\Vert \geq e^{\bar CN} \qquad(\bar C>0)\tag 5.2
$$
and there is $\ell \in\Bbb Z^{d+\nu}, |\ell|\sim N^\tau \,(\tau>7C(d+\nu)\nu)$,
$ \Lambda_N(\ell)$ satifying
$$
(\Lambda_N(\ell) \cap \Bbb Z^d)\cap (\Lambda_{\bar N}\cap \Bbb Z^d)=\emptyset\tag 5.3
$$
such that
$$
|(H_{\Lambda_N(\ell)}-E)^{-1} (m, m')|>e^{-\gamma|m-m'|}\tag 5.4
$$
for some $m, m'\in \Lambda_N(\ell), |m-m'|> N/4$.
$$
\text{\rm  Prob} \leq N^{-q} \qquad \big(q>p-\tau(d+\nu)-1\big)
\tag 5.5
$$
\endproclaim

\noindent
{\bf Proof.} (5.3) implies that $H_{\Lambda_N(\ell)}$ is independent from $H_{\bar N}$ and hence its eigenvalues.
Summing over the probabilistic estimates in (3.5), we obtain (5.5).
(The lemma holds as soon as (5.3) is satisfied due to independence.
We take $|\ell|\sim N^\tau$ in view of Lemma 4.3).
$\hfill\square$

\noindent
{\bf Proof of the theorem.}

Let $X_\epsilon \subset \Bbb R^{\Bbb Z^d}$ be the probability subspace defined in (4.4), i.e., Proposition 2.7 with the
modification 
$ i\in [-\bar N, \bar N]^d, \bar N\asymp N^C$, (2.63) still holds in this case.

Let $\tilde X_N$ be the complement of the set defined in (5.1).
Let
$$
\tilde X_{N_0} =\big(\bigcap_{N>N_0} \tilde X_N \big)\bigcap \tilde X_\epsilon.
$$

Fix $x\in \tilde X_{N_0}$.
Let $\Omega_{N, x}$ be the complement of the set defined in (4.26).
Let
$$
\Omega_{x, N_0} =\bigcap_{N>N_0} \Omega_{N, x}.
$$
For any given $\eta>0$, $\exists\, N_0$, such that
$$
\tilde X_{\eta, \ve} \overset\text{def}\to= \tilde X_{N_0}.
$$
Satisfying
$$
\mes \tilde X_{\eta, \ve} \geq 1-\eta
$$
and
$$
\Omega_{x, \eta, \ve}\overset\text{def}\to= \Omega_{x, N_0}.
$$
Satisfying
$$
\mes \Omega_{x, \eta, \ve}\geq 1-\eta,
$$
if $p'> 422(d+\nu)\nu $ (cf. A1, A2, 2.69, 3.12, 3.18, 4.25, 4.28, 5.5).

The proof uses lemmas 4.3 and 5.1 and follows the same strategy as in \cite{BG, BGS}.
So we will only highlight the main points.

\itemitem {$\bullet$} \ Generalized eigenfunctions of $H:H\psi =E\psi$ has the apriori bound
$$
|\psi(m)|\leq 1+|m|^c\tag 5.6
$$
from the Schnol-Simon Theorem [Sh, Sim].

\itemitem {$\bullet$} \ Let $\psi$ be a non-zero eigenfunction of $H: H\psi = E\psi$.
Let $\Lambda\subset\Bbb Z^{d+\nu}$ and assume $E\not\in \sigma (H_\Lambda)$.
Then for all $m\in \Bbb Z^{d+\nu}$
$$
\psi(m) =\sum_{\{m'\in \Lambda|\exists \, m''\in \Bbb Z^{d+\nu}\backslash\Lambda, |m'-m''|=1\}}
(H_\Lambda-E)^{-1} (m, m')\psi (m'').\tag 5.7
$$

\itemitem{$\bullet$}
 \ From semi-algebraic considerations and the restriction of $\omega\in\Omega_{x, \eta, \epsilon} $ (see Lemma 2.3),
$\exists \, \kappa>0, M\asymp \bar N$ $(\bar N\asymp N^C, C>10(d+\nu)$, cf. Lemma 2.6) such that
$$
\align
&\#\{j\in\Bbb Z^{d+\nu}|\Lambda_N(j)\cap \Lambda_M\not=\emptyset, \Lambda_N(j) \text { is a bad $N$-box}\}
\\
&\leq M^{1-\kappa}.\tag 5.8
\endalign
$$
Let $\Cal I$ be the set defined in (5.8).
For $C$ large enough, there exists an annulus
$$
A=\Lambda_L\backslash \Lambda_{L'}\subset \Lambda_M
$$
of width $L'-L>2N$ such that
$$
A\cap \Cal I=\emptyset.
$$

\itemitem{$\bullet$}\ Without loss, assume $\psi(0)=1$.
Using (5.6, 5.7) first with $\Lambda=\Lambda_N(j), j\in A$ and then $\Lambda=\Lambda_{\bar N}$ we obtain
$$
\Vert(H_{\bar N}-E)^{-1}\Vert\geq e^{\bar CN}\qquad (\bar C>0).\tag 5.9
$$

\itemitem{$\bullet$}\
Lemma 4.3, 5.1 and an application of the resolvent equation as used earlier imply that $\forall i \in\partial
\Lambda_{2K}(0)$, $K\sim N^{\tau}$
$$
|G_{\Lambda_K(i)} (E, m, m')|\leq \exp (-|m-m'|).\tag 5.10
$$
$m, m'\in \Lambda_K(i), |m-m'|\geq K/4$.

We note that
$$
\bigcup_{i\in \partial \Lambda_{2K}(0)} \Lambda_K(i) =\Lambda_{3K}(0)\backslash \Lambda_K(0)\overset \text{def}\to=
U
$$
(5.6, 5.7, 5.10) imply that
$$
|\psi(m)| < e^{-|m|/2}
$$
for $m$ such that
$\dist (m, \partial U)\geq K/4$, provided $N$ and thus $K$ are large.
$\hfill\square$

\heading
6. Proof of Anderson localization for wave operator
\endheading

The quasi-energy operator $K_w$ in the wave case is
$$
K_w =-\sum^\nu_{k=1} \sum^\nu_{ k'=1} \omega_k\omega_{k'} \, \frac {\partial^2}{\partial\theta_k\partial\theta_{k'}}
+\epsilon\Delta+ V+ \sum^\nu_{k=1} W_k \cos 2\pi \theta_k\tag 6.1
$$
on $\ell^2(\Bbb Z^d)\times L^2(\Bbb T^\nu)$, where $\omega =(\omega_1, \omega_2\cdots \omega_\nu)\in (0, 1]^\nu$. 
$V$ is the random potential on $\Bbb Z^d$, $0<\epsilon \ll 1$, and $W_k$ satisfies the decay properties specified in 
(H1).
Compared with the quasi-energy operator $K$ for Schr\"odinger in (2.1), the only difference is that the $\theta$
deriviatives are second order.

Performing a partial Fourier series transform in the $\Bbb T^\nu$ variables as in sect. 2, we are led to study the
following unitarily equivalent operator:
$$
H_w =\delta_j\tilde\Delta_n+(n\cdot\omega)^2 +\epsilon \Delta_j+V_j\tag 6.2
$$
on $\ell^2(\Bbb Z^{d+\nu})$, where $n\in \Bbb Z^\nu, j\in\Bbb Z^d$ and $\delta_j\tilde\Delta_n$ is as in (2.3).
We proceed as in the Schr\"odinger case and inroduce the parameter $\theta \in\Bbb R$.
We define
$$
H_w(\theta)=\delta_j\tilde\Delta_n+ (n\cdot \omega+\theta)^2+\epsilon\Delta_j+V_j\tag 6.3
$$
on $\ell^2(\Bbb Z^{d+\nu})$ and study the Green's functions
$$
G_{w, \Lambda}(\theta, E)=(H_{w,\Lambda}(\theta)-E)^{-1}\tag 6.4
$$
for a class of finite sets $\Lambda\subset\Bbb Z^{d+\nu}$ $\Lambda\nearrow \Bbb Z^{d+\nu}$ to be specified shortly.

The main difference between $H_w(\theta)$ in (6.3) and $H(\theta)$ in (2.5) is that
$$
\frac{\partial H(\theta)}{\partial\theta}=1; \tag 6.5
$$
while
$$
\frac{\partial H_w(\theta)}{\partial \theta}=2(n\cdot\omega+\theta)\tag 6.6
$$
which could be 0.
So the apriori estimate in Lemma 2.2 for Schr\"odinger does not apply here.
($\theta$ and $E$ are no longer equivalent, see (2.17).)
We need to resort to Cartan type of theorem for analytic matrix-valued functions as in \cite{BGS}.
Unlike Lemma 2.2, which holds at all scales this requires a multi-scale analysis.
At each scale $N(N\gg 1)$ we need measure estimates on the bad sets at two previous scales $N_0, N_1$ with
$$
\log N_0<\log N_1\ll \log N.
$$

We extend the class of finite subsets $\Lambda$ of $\Bbb Z^{d+\nu}$, which were previously cubes to {\it elementary
regions}: (as in \cite{BGS})
$$
\Lambda\overset\text{def}\to= R\backslash ( R+m)\tag 6.7
$$
where $m\in \Bbb Z^{d+\nu}$ is arbitrary and $ R$ is a rectangle
$$
 R=\times^{d+\nu}_{i=1} [-M_i, M_i]+k,\tag 6.8
$$
$k\in\Bbb Z^{d+\nu}$. 
The {\it size} of $\Lambda$, denoted by $\ell(\Lambda)$, is simply its diameter.
We denote by $\Cal {ER}(M)$, the set of all elementary regions of size $M$.

Let
$$
\sigma\in(0, 1), N=[N_0^C]+1 \qquad (C>1)\tag 6.9
$$
$N_0\gg 1$ (determined by $\delta, \sigma$ similar to (2.11), $\delta$ as in (2.31),
$$
\Lambda\subset [-N, N]^{d+\nu}, \Lambda\in\Cal {ER}(N).
$$
Let $\Lambda, \Lambda'$ be two elementary regions.
Let $\tilde\Lambda, \tilde{\Lambda'}$ be their respective convex envelop.
We say that $\Lambda$ and $\Lambda'$ are {\it disjoint} if $\tilde \Lambda\cap \tilde{\Lambda'} =\emptyset$.
Following is the analogue of Lemma 2.3.

\proclaim
{Lemma 6.1} 
Fix $x\in X_N$, $X_N$ as in Sect. 2.
There exists a set $\Omega_N\subset (0, 1]^\nu$,
$$
\text{\rm mes\,}\Omega_N\geq 1-e^{N^{\frac \sigma {2C}}}
$$
such that if $\omega\in\Omega_N$, then for {\it any}
fixed $\theta, E$, there are at most {\it two} pair-wise disjoint bad $\Lambda_0$'s,
$$
\Lambda_0\in \bigcup_{N_0\leq M\leq 2N_0} \Cal {ER}(M)
$$
in $\Lambda$.
Moreovr, $(0, 1]^\nu\backslash \Omega_N$ is contained in the unon of at most
$N^{5(d+\nu)\nu}$ components.
\endproclaim

\noindent
{\bf Proof.}
We follow the same line of argument as in the proof of
Lemma 2.3.
Let
$$
\Lambda_0, \Lambda_0', \Lambda_0''\in \bigcup_{N_0\leq M\leq 2N_0} \Cal {ER}(M)
$$
$\Lambda_0, \Lambda_0', \Lambda_0''\subset\Lambda$, be pair-wise disjoint.

Let 
$$
\text{$\Lambda_{0, j}$ be the projection of $\Lambda_0$ onto $\Bbb Z^d$}.
$$
$$
\Lambda_{0, n}  \text { be the projection of $\Lambda_0$ {onto} $\Bbb Z^\nu$}
$$
and similarly for $\Lambda_0', \Lambda_0''$.

Assume $\Lambda_0, \Lambda_0'$ and $\Lambda_0''$ are all bad, then ther exist
$$
\matrix
& n\in\Lambda_{0, n}   &\qquad n'\in \Lambda_{0, n}'  &\qquad n''\in \Lambda_{0, n}''\\ 
&j\in \Lambda _{0, j} &\qquad j'\in \Lambda'_{0, j}     & \qquad j''\in \Lambda_{0, j}''
\endmatrix
\qquad \text { such that }
$$
$$
\cases
|(n\cdot\omega+\theta)^2 +\mu_j-E|&\leq 2e^{-N_0^\sigma}\\
|(n'\cdot\omega+\theta)^2+\mu_{j'} -E|&\leq 2e^{-N_0^\sigma}\\
|(n''\cdot\omega+\theta)^2+\mu_{j''}-E|&\leq 2e^{N_0^\sigma}\endcases
\tag 6.10
$$
Assume $n, n', n''$ are distinct, otherwise say $n=n'$, then the first two inequalities in (6.10) imply that
$$
(\mu_j-\mu_j)< 4e^{-N_0^\sigma}\tag 6.11
$$
where $\mu_j$ is an eigenvalue for $\Lambda_{0, j}$ and $\mu_j', \Lambda_{0, j}'$.
Since $\Lambda_{0, j}\cap \Lambda_{0, j'}=\emptyset$, Theorem A, in particular (A3) implies that $|\mu_j-\mu_j'|\geq
e^{-N_0^\beta}\gg e^{-N_0^\sigma}$ by choosing $\sigma>\beta$ as in the proof of Lemma 2.3.

Subtracting the inequalities in (6.10) pairwise, we get two inequalities with linear dependence on $\theta$ and
independent of $E$.
Eliminating the dependence on $\theta$, we obtain the following:
$$
\align
&|(n-n') \cdot\omega(n-n'')\omega(n'-n'')\cdot \omega+ (n-n'')\cdot \omega (\mu_j-\mu_j')\\
&-(n-n')\cdot\omega(\mu_j-\mu_j;'')| \leq 4 e^{-N_0^\sigma} \big(|(n-n')\cdot\omega|+
|(n-n'')\cdot \omega|).\tag 6.12
\endalign
$$
Since $\omega\in (0, 1]^\nu, n-n', n-n''\in [-2N, 2N]^\nu$
$$
\align
\text {RHS of}\, (6.12) &\leq 16N\nu E e^{-N^\sigma_0}\\
&\leq e^{-\frac{N_0^\sigma}{2}}\tag 6.13
\endalign
$$
for  appropriate $\sigma, C$.
(Recall $N=[N_0^C]+1)$.

Let 
$$
\cases
m=n-n', m'=n-n''\\
\lambda =\mu_j-\mu_j', \lambda'= \mu_j-\mu_j''.\endcases\tag 6.14
$$
The solutions to (6.12) is contained in the solutions to
$$
|(m\cdot\omega)(m'\cdot\omega)(m-m')\cdot\omega+(\lambda m'-\lambda' m)\omega|\leq e^{-\frac{N_0^\sigma}2}
\quad \ (m\not= m', m\not=0, m'\not=0).\tag 6.15
$$
Assume $\omega$ is a solution to
$$
(m\cdot\omega)(m'\cdot\omega)(m-m')\cdot\omega+(\lambda m'-\lambda' m)\cdot\omega=0 \
\quad (m\not=m', m\not=0, m'\not=0)\tag 6.16
$$

Let $\omega\to \omega+\delta (|\delta|\ll 1)$.
It is easy to see that we can always choose $\delta$ so that the third order variation
$$
|(m\cdot\delta)(m'\cdot\delta)(m-m')\cdot\delta|> c'|\delta|^3\tag 6.17
$$
for all $m\not=0, m'\not= 0, m'\not=m'$, where $c'>0$ is independent of $m, m'$.
We obtain
$$
\align
&\mes \{\omega \in (0, 1]^\nu| (6.15) \text { is satisfied}\}\\
&\leq Ce^{-\frac {N_0^\sigma}6}\tag 6.18
\endalign
$$ where $C$ only depends on $\nu$.

There are at most $N^{2d+5\nu}$ equations of the form (6.15).
Let
$$
\Omega_N=\{\omega\in (0, 1]^\nu| (6.15) \text { is satisfied } \forall m, m', \lambda,\lambda'\}
$$
We obtain the lemma by using (6.18) and Basu's theorem stated as Theorem 7.3 in \cite{BGS}.
$\hfill\square$

Assume $N\gg 1, \Lambda\in \Cal E R(N), \Lambda\subset [-N, N]^{d+\nu}$.
Let $X\subset\Bbb R^{\Bbb Z^d}$ be defined similarly as in Lemma 2.6.
For any $x\in X$, define $\Omega_x$ similarly to that in Lemma 2.6.
Combining Lemma 6.1 and semi-algebraic considerations, we have as in Lemma 2.6 that there exists $\delta_0>0$, such
that for any fixed $\theta, x\in X, \omega\in\Omega_x\bigcap DC_{A, c}(2N)$
$$
\align
\# \{m\in[-N, N]^{d+\nu}|\exists\, \Lambda_0 \in \Cal ER(M), N_0&\leq M\leq 2N_0,
\Lambda_0\subset m+[-M,M]^{d+\nu}, \Lambda_0 \text { is bad}\}\\
&\leq N^{1-\delta_0}.\tag 6.19
\endalign
$$
where $N_0=[N^{1/C}]$ for appropriate $C$ depending on $\delta_0$ only.

Recall the definition of good and bad regions for fixed $\theta$.
$\Lambda_0$ of size $\ell(\Lambda_0)$ is good if
$$
\align
&\Vert G_{\Lambda_0}(\theta, E)\Vert< e^{N_0^\sigma}\\
&|G_{\Lambda_0}(\theta, E)(m, m')| < e^{-\gamma|m-m'|}\tag 6.20
\endalign
$$
for all $m, m'\in \Lambda_0$, $|m-m'|>\ell (\Lambda_0)/4$ where $\sigma>0, \gamma>0$.
Otherwise it is bad. As in (2.38), $\Cal G_x^{\gamma, \sigma} (\Lambda_0, E)$ is the set of $\theta\in\Bbb R$
such that (6.20) holds and $\Cal B_x^{\gamma, \sigma} (\Lambda_0, E)$ is the complement set.

The following lemma plays the role of Lemma 2.2 for $H_w$.

\proclaim
{Lemma 6.2} Assume $\epsilon,\delta\ll 1$ and (2.3).
There exist $\sigma, \rho, \gamma>0$ satisfying $0<\sigma, \rho < 1, \sigma+\delta_0> 1+3\rho$, where $\delta_0$ is as
in (6.19), and $C_1(\sigma, \rho)\gg\frac 1\rho, N_0\leq N_1$, satisfying
$$
\bar N_0(\gamma, \sigma,\rho)\leq 100 N_0\leq N_1^\rho,\tag 6.21
$$
such that for any $N_0\leq M\leq N$, and any $\Lambda\in \Cal ER(M)$
$$
\sup_{x\in X, E} \text{\rm mes\,} \big(\Cal B_\theta^{\gamma, \sigma}(\Lambda, E)\big)
\leq \exp (-\ell (n)^\rho).\tag 6.22
$$
Assume moreover that
$$
\omega\in\Omega_x\cap DC_{A, c}(2N_1^{\rho C_1}).
$$
Then for all $\Lambda\subset\Cal {ER}(N)$
$$
\sup_{x\in X, E} \text{\rm mes\,} \{\theta\in\Bbb R\big| \, \Vert G_\Lambda (\theta, E)\Vert> e^{N^b}\}< e^{-N^{3\rho}}
\tag 6.23
$$
where $N_0^{C_1} \leq N\leq N_1^{\rho C_1}$.
\endproclaim

The proof of the above lemma is very similar to the proof of Lemma 4.4 of \cite{BGS}, (see also Chap XIV of \cite{Bo}
for a more detailed exposition). So instead of replicating the proof, we only sketch the main line of arguments.

\itemitem {$\bullet$} As mentioned earlier in (6.5, 6.6), contrary to the Schr\"odinger case, the first order
variation can vanish. So we need to resort to analytic and subharmonic function theory to control the
measure of $\Cal B_x^{\gamma, \sigma} (\Lambda, E)$. To do that we need 2 scales $M$, $M_1$, with 
$\log M< \log M_1\ll \log N$.

\itemitem {$\bullet$} Fix $\theta$, at scale $M$ ($N_0\leq M\leq 2N_0$), let $\Lambda_{*}$ be, roughly 
speaking, the complement of the set in (6.19). For more precise definition, which requires a partition
of $\Lambda$, see the beginning of the proof of Lemma 4.4 in \cite{BGS}. Using an elementary resolvent
expansion (Lemma 2.2 of \cite{BGS}), we obtain an upper bound on $\Vert G_{\Lambda_{*}}(E,\theta)\Vert$
by using the decay estimate on the $\Lambda_0$'s, elementary regions at scale M, in $\Lambda_{*}$. By 
definition they are all {\it good}. By standard Neumann series arguments, this bound is preserved inside
the disk $B(\theta,e^{-N_0})\subset\Bbb C$. 

\noindent{\it Remark.} We have control over the size of $\Lambda_{*}$ via (6.19), but not its geometry. Typically
$\Lambda_{*}$ is {\it non-convex}. Hence the need for elementary regions which are more general than cubes, in
particular L-shaped regions, in view of Lemma 2.2 of \cite{BGS}. 

\itemitem {$\bullet$} Define a matrix-valued analytic function $A(\theta')$ on $B(\theta, e^{-N_0})$ as
$$\aligned A(\theta')&=R_{\Lambda_{*}^c} H_w(\theta') R_{\Lambda_{*}^c}\\
&\quad -R_{\Lambda_{*}^c} H_w(\theta') R_{\Lambda_{*}}G_{\Lambda_{*}}(E,\theta')
R_{\Lambda_{*}} H_w(\theta') R_{\Lambda_{*}^c}\endaligned\tag 6.24$$
where $\Lambda_{*}^c=\Lambda\backslash \Lambda_{*}$, $R_{\Lambda_{*}}$, $R_{\Lambda_{*}}^c$ are projections.
From (6.19), $A(\theta')$ is a rank $\Cal O(N^{1-\delta_0})\times \Cal O(N^{1-\delta_0})$ matrix. The
raison d'etre of introducing $A(\theta')$ is the following inequality:
$$\Vert A(\theta')^{-1}\Vert\lesssim \Vert G_\Lambda(\theta',E)\Vert \lesssim e^{2N_0}\Vert A(\theta')^{-1}\Vert,\tag
6.25$$ (see Lemma 4.8 of \cite{BGS}). So to bound $\Vert G_\Lambda(\theta',E)\Vert $, it is sufficient to bound
$\Vert A(\theta')^{-1}\Vert$, which is of smaller dimension.

\itemitem {$\bullet$} Toward that end, we introduce an intermediate scale $M_1$, $\log M_1>\log M$. We work in an 
interval $\Theta=\{\theta'||\theta'-\theta|<e^{-N_0}\}$. Using (6.21) for the $\Lambda_1$'s at scale $M_1$ and in
$\Lambda$, the same elementary resolvent expansion, we obtain a bound on $\Vert G_\Lambda(\theta',E)\Vert$ except
for a set of $\theta'$ of measure smaller than $e^{-\Cal O(N_1^\rho)}$. So there exists $y\in\Theta$, such
that we have {\it both} a lower bound on $\Vert A(\theta')\Vert$ at $\theta'=y$, and an apriori upper bound on
$\Vert A(\theta')\Vert$ in the disk $B(y,e^{-N_0}/2)$, which comes from boundedness of $H_w$ and the bound on
$\Vert G_{\Lambda^{\*}}(E,\theta)\Vert$ (see (6.24)).

\itemitem {$\bullet$} Transfering the estimates on $\Vert A(\theta')\Vert$ into estimates on $\log|\det A(\theta')|$,
which is subharmonic and using either Cartan type of theorem (see sect. 11.2 in \cite {Le}) or proceeding as
in the proof of Lemma 4.4 of \cite{BGS} or Chap XIV of \cite{Bo}, we obatin the lemma by covering the interval
$I=(-\Cal O(N_1^{\rho C_1},\Cal O(N_1^{\rho C_1}))$ with intervals of size $e^{-N_0}$. (Recall (6.21) and
that for all $\theta\notin I$, $H_{w,\Lambda}-E$ is automatically invertible.)
$\hfill\square$

Lemma 3.1 remains valid for $H_{w,\Lambda}$, as the first order variation in $x$ remains to be $1$.
We can now proceed as in the Schr\"odinger case to prove Anderson localization for the wave operator $H_w$.
We obtain the same Theorem as in sect. 5 with $H_w$ in place of $H$.

\bigskip
\heading
Appendix:
Localization results for Random Schr\"odinger operators
\endheading

Random Schr\"odinger opeator is the operator
$$
H=\epsilon \Delta +V \text { on } \ell^2 (\Bbb Z^d)
$$
where $\epsilon>0$ is a parameter, $\Delta (i, j) =1$ if $|i-j|=1$ and zero otherwise, $V=\{v_i\} _{i\in\Bbb Z^d}$ is a
family of independent identically distributed (iid) random variables with common probability distribution $g$.
The spectrum of $H$ is given by
$$
\align
\sigma (H) &=\sigma(\epsilon\Delta)+\sigma(V)\\
&=[-2\epsilon d, 2\epsilon d]+ \supp g.
\endalign
$$

There are a few versions of Anderson localization results for $H$,
The one that is most adapted for our purpose is proven in \cite{vDK}, which we restate below.

For any $L\in\Bbb N$, let $\Lambda_L(i)= [-L, L]^d+i, i\in\Bbb Z^d$.
Let $m>0, E\in \Bbb R$.
$\Lambda_L(i)$ is $(m, E)$-regular (for a fixed $V$) if $E\not\in \sigma (H_{\Lambda_l(i)})$ and
$$
|G_{\Lambda_L(i)}(E; j, j')|\leq e^{-m|j-j'|}
$$
for all $j\cdot j'\in \Lambda_L(i), |j-j'|> L/4$.

\proclaim
{Theorem A}
Let $I\subset\Bbb R$ be a bounded interval.
Suppose that for some $L_0>0$, we have
$$
 \text {{\rm Prob} \{for any $E\in I$ either $\Lambda_{L_0}(i)$ or $\Lambda_{L_0}(j)$ is
$(m_0, E)$-regular\} $\geq 1-\frac 1{L_0^{2p'}}$},\tag A1
$$ 
for some $p'>d, m_0>0$, and any $i, j\in\Bbb Z^d,  |i-j|> 2L_0$
$$
 \text{\rm Prob\,}\{\text{\rm dist\,} \big(E, \sigma(H_{\Lambda_L(0)})\big) < e^{-L^\beta}\}\leq 1/{L^q}\tag A2
$$
for some $\beta$ and $q$, $0<\beta<1$ and $q> 4p+6d$ all $E$ with
$$
\text{\rm dist\,} (E, I)\leq \frac 12 e^{-L^\beta},
$$
and all $L\geq L_0$.
Then there exists $\alpha, 1<\alpha <2$, such that if we set $L_{k+1}=[L^\alpha_k]+1$, $k= 0, 1,2 \ldots$ and pick $m,
0<m<m_0$, there is $Q<\infty$, such that if $L_0>Q$, we have that for any $k=0, 1,2 \ldots$
$$
\text{{\rm Prob} \{for any $E\in I$ either $\Lambda_{L_k(i)} $ or $\Lambda_{L_k(j)}$ is $(m, E)$ regular\} 
$\geq 1-\frac 1{L_k^{2p}}$}
$$
for any $i, j\in\Bbb Z^d$ with $|i-j|> 2L_k$.
\endproclaim

\noindent
{\it Remark.}
On the same probability subspace,
$$
\dist \big(\sigma\big(H_{\Lambda_{L_k}(i)}\big), \sigma \big(H_{\Lambda_{L_k}(j)}\big)\big)> e^{-L_k^\beta},\qquad
\beta>0\tag A3
$$
if $|i-j|> 2L_k$.
This is part of the ingredient of the proof of Theorem A.

(A1) is verified if $\epsilon$ is sufficiently small. (A2) is provided by the Wegner Lemma if $g$ is absolutely
continuous with a bounded density $\tilde g$.
$$
\{\text {Prob}\big(\dist \big(E, \sigma(H_\Lambda)\big)\leq K\big)\} \leq CK|\Lambda|\Vert\tilde g\Vert_\infty.
$$

\end